\documentclass[10pt]{article}
\usepackage[utf8]{inputenc}
\usepackage{epsfig}
\usepackage{color}
\usepackage[british,english]{babel}
\usepackage{amsthm}
\usepackage{amsmath}
\usepackage{amsfonts}
\usepackage{amssymb}
\usepackage{graphicx}
\newtheorem{corollary}{Corollary}[section]
\newtheorem{definition}[corollary]{Definition}

\newtheorem{lemma}[corollary]{Lemma}
\newtheorem{proposition}[corollary]{Proposition}
\newtheorem{remark}[corollary]{Remark}
\newtheorem{theorem}[corollary]{Theorem}
\newfont{\sBlackboard}{msbm10 scaled 900}

\newcommand{\mylabel}[1]{\label{#1}
            \ifx\undefined\stillediting
            \else \fbox{$#1$}\fi }
\newcommand{\BE}{\begin{equation}}

\newcommand{\EEQ}{\end{equation}}
\newcommand{\rfb}[1]{\mbox{\rm
   (\ref{#1})}\ifx\undefined\stillediting\else:\fbox{$#1$}\fi}

\newfont{\Blackboard}{msbm10 scaled 1200}

\newfont{\roma}{cmr10 scaled 1200}

\def\CC{\rm \hbox{C\kern-.56em\raise.4ex
         \hbox{$\scriptscriptstyle |$}\kern+0.5 em }}

\def\n{|\kern -.05cm{|}\kern -.05cm{|}}

\def \noame{\noalign{\medskip}}
\newcommand{\mm}    {{\hbox{\hskip 0.5pt}}}

\newcommand{\bluff} {{\hbox{\raise 15pt \hbox{\mm}}}}

\newcommand{\ep}   {\varepsilon}
\usepackage{fancyhdr}

\makeatletter
\def\section{\@startsection {section}{1}{\z@}{-3.5ex plus -1ex minus
    -.2ex}{2.3ex plus .2ex}{\large\bf}}
\makeatother
\def\be{\begin{equation}}
\def\ee{\end{equation}}

\date{ }
\begin{document}
\thispagestyle{empty}
\title{\Large \bf  Two-dimensional Carreau law for a quasi-newtonian fluid  flow through a thin domain with a slightly rough boundary}\maketitle
\vspace{-2cm}
\begin{center}
Mar\'ia ANGUIANO\footnote{Departamento de An\'alisis Matem\'atico. Facultad de Matem\'aticas. Universidad de Sevilla. 41012-Sevilla (Spain) anguiano@us.es} and Francisco Javier SU\'AREZ-GRAU\footnote{Departamento de Ecuaciones Diferenciales y An\'alisis Num\'erico. Facultad de Matem\'aticas. Universidad de Sevilla. 41012-Sevilla (Spain) fjsgrau@us.es}
 \end{center}

 \renewcommand{\abstractname} {\bf Abstract}
\begin{abstract}
This study investigates the asymptotic behavior of the steady-state quasi-Newtonian Stokes
flow with viscosity given by the Carreau law within a thin domain, focusing on the effects of a  rough boundary of the domain.
Employing asymptotic techniques with respect to the domain's thickness,  we rigorously derive the effective nonlinear two-dimensional Reynolds 
model describing the fluid flow. The mathematical analysis is based on deriving the sharp a priori
estimates and proving the compactness results of the rescaled functions together with monotonicity arguments. The resulting limit model incorporates contributions of the oscillating boundary and thus, it could prove useful in the applications involving this lubrication regime.
\end{abstract}
\bigskip\noindent

\noindent {\small \bf AMS classification numbers:}  35B27; 35Q35; 76A05; 76M50; 76A20.  \\

\ \\
\section {Introduction}\label{S1}
The classical lubrication problem is to describe the situation in which two adjacent surfaces in relative motion are separated by a thin film of fluid acting as a lubricant. Such situation appears naturally in
numerous industrial and engineering applications. The mathematical models for describing the motion of the lubricant usually result from the
simplification of the geometry of the lubricant film, i.e. its thickness. Using the film thickness as a small
parameter, an asymptotic approximation of the Stokes system can be derived providing the well-known Reynolds equation for the pressure of the fluid, see Bayada \& Chambat \cite{Bayada_Chambat_Transition} for more details. Thus, in the stationary case, and considering no-slip condition on the top and bottom boundaries and an exterior force ${\bf f'}:x'\in\mathbb{R}^2 \to \mathbb{R}^2$, in the Newtonian case (where the  viscosity is constant), the two-dimensional Reynolds equation for the pressure $\widetilde p$ has the following form 
\begin{equation}\label{ReynoldsNew}
{\rm div}_{x'}\left({h(x')^3\over 12\mu}\left({\bf f'}(x')-\nabla_{x'}\widetilde p(x')\right)\right)=0,
\end{equation}
where $h$ describes the shape of the top boundary and $\mu$ is the viscosity of the fluid.

 However, many others fluids behave differently, the viscosity of these fluids is no more constant, such as pastes or polymer solutions. These fluids  are called non-Newtonian (see Ganji \cite{Ganji} for more details).  Non-Newtonian fluids are common in polymer processing, biological fluids, paints and coatings, food industry, and lubricants. Modeling their flow in confined/thin geometries is critical for predicting performance of devices such as bearings, medical devices, extrusion dies, and thin films, see for instance Agassant {\it et al.} \cite{Agassant}, Barnes {\it et al.} \cite{Barnes} or Bird {\it et al.} \cite{Bird}. Mathematically, nonlinear constitutive laws complicate existence, uniqueness, and asymptotic analysis (see for instance Mikeli\'c \cite{Mikelic2}. Understanding non-Newtonian fluids is vital in environmental engineering,  so  homogenization and asymptotic expansion methods are essential to derive reduced models for thin domains or porous media (see for instance Anguiano \& Su\'arez-Grau \cite{Anguiano_SG},  Boughanim \& R. Tapi\'ero \cite{Tapiero2}, Bourgeat \& Mikeli\'c  \cite{Bourgeat1} or Bourgeat {\it et al.} \cite{Bourgeat2}).

The simplest idea to describe non-Newtonian fluids is to plot the viscosity measurements versus the imposed shear rate and then, to fit the obtained curve with a simple template viscosity function, adjusting some few parameters. This is the main idea of quasi-Newtonian fluids models, which could be viewed as a first step inside the world of non-Newtonian fluids models (see Saramito \cite[Chapter 2]{Saramito} for more details).

The incompressible quasi-Newtonian fluids are characterized by the viscosity depending on
the principal invariants of the symmetric stretching tensor $\mathbb{D}[{\bf u}]$. If ${\bf u}$ is the velocity, $p$ is the pressure and $D{\bf u}$ is the gradient velocity tensor, $\mathbb{D}[{\bf u}]=(D{\bf u}+D^t {\bf u})/2$ denotes the symmetric stretching tensor and $\sigma$ the stress tensor given by $\sigma=-pI+2\eta_r \mathbb{D}[{\bf u}]$. The viscosity $\eta_r$ is constant for a Newtonian fluid but dependent on the shear rate, that is, $\eta_r=\eta_r(\mathbb{D}[{\bf u}])$, for viscous non-Newtonian fluids. The deviatoric stress tensor $\tau$, that is, the part of the total stress tensor that is zero at equilibrium, is then a nonlinear function of the shear rate $\mathbb{D}[{\bf u}]$, that is $\tau=\eta_r(\mathbb{D}[{\bf u}])\mathbb{D}[{\bf u}]$, see Barnes {\it et al.} \cite{Barnes}, Bird {\it et al.} \cite{Bird} and Mikeli\'c \cite{MikelicIntro1,MikelicIntro2} for more details).

A widely used law in engineering practice is the power law model  (see for instance Agassant {\it et al.} \cite{Agassant}, Barnes {\it et al.} \cite{Barnes},  Bird {\it et al.} \cite{Bird}) or Chhabra \& Richardson \cite{Chhabra}),  where the viscosity as a function of the shear rate is given by 
$$\eta_r(\mathbb{D}[{\bf u}])=\mu |\mathbb{D}[{\bf u}]|^{r-2},\quad 1<r<+\infty,$$
where the two material parameters $\mu>0$ and $r$ are called the consistency and the flow index, respectively. Here, the matrix norm $|\cdot |$ is defined by $|\xi|^2=Tr(\xi\xi^t)$ with $\xi\in \mathbb{R}^3$.  We recall the classification in terms of $r$:
\begin{itemize}
\item[--] Shear-thinning (Pseudoplastic) ($1<r<2$): viscosity decreases with increased shear rate, e.g., blood, polymer solutions.
\item[--] Newtonian fluids ($r=2$): constant viscosity.
\item[--] Shear-thickening (Dilatant) $r>2$: viscosity increases with more shear, e.g., some suspensions.
 
\end{itemize}
Similarly to the derivation of equation (\ref{ReynoldsNew}), a two-dimensional nonlinear Reynolds equation for quasi-Newtonian fluid with viscosity given by the power law has been obtained in Mikeli\'c \& Tapiero \cite{Tapiero}, with the following form
\begin{equation}\label{ReynoldsPowerLaw}
{\rm div}_{x'}\left(
{h(x')^{r'+1}\over 2^{r'\over 2}(r'+1)\mu^{r'-1}}\left|{\bf f}'(x')-\nabla_{x'}\widetilde p(x')\right|^{r'-2}\left({\bf f}'(x')-\nabla_{x'}\widetilde p(x')\right)\right)=0,
\end{equation}
where $r'$ is the conjugate exponent of $r$, i.e. $1/r+1/r'=1$.

It is observed that the power law correctly describes the behavior of polymers at high-shear rates. It offers the advantage
of allowing analytical calculations in simple geometries. However, it has the disadvantage of not describing a Newtonian plateau and even predicts an infinite viscosity as the shear rate goes to zero
and $1<r<2$  (see Agassant {\it et al.} \cite[p. 49]{Agassant} for more details), whereas for real fluids, it tends to some constant
value $\eta_0$ called the zero-shear-rate viscosity. Thus, to solve this problem, it is considered  the Carreau law,  which is given by
\begin{equation}\label{Carreaulaw}\eta_r(\xi)=(\eta_0-\eta_\infty)(1+\lambda|\xi|^2)^{{r\over 2}-1}+\eta_\infty,\quad 1<r<+\infty,\quad\eta_0>\eta_\infty>0,\quad \lambda>0,
\end{equation}
where $\eta_\infty$ is the high-shear-rate limit of the viscosity, the parameter $\lambda$ is a time constant and $r-1$ is a dimensionless constant describing the slope in the power law region (see Saramito \cite[Chapter 2]{Saramito} for more details). Also,  a two-dimensional nonlinear Reynolds equation for quasi-Newtonian fluid with viscosity given by the Carreau law has been obtained in Boughanim \& Tapiero \cite{Tapiero2}, with the following form
\begin{equation}\label{ReynoldsCarreauLaw}
{\rm div}_{x'}\left\{\left(\int_{{-h(x')\over 2}}^{h(x')\over 2}{\left({h(x')\over 2}+\xi\right)\xi\over \psi(2|{\bf f}'(x')-\nabla_{x'}\widetilde p(x')||\xi|)}d\xi\right)\left({\bf f}'(x')-\nabla_{x'}\widetilde p(x')\right)
\right\}=0,
\end{equation}
where $\psi$ is the inverse function of the algebraic equation
\begin{equation}\label{algebraic}
\tau=\zeta\sqrt{{2\over \lambda}\left\{{\zeta-\eta_\infty\over \eta_0-\eta_\infty}\right\}-1},
\end{equation}
which has a unique solution noted $\zeta=\psi(\eta)$ for $\eta\in\mathbb{R}^+$, see \cite[Proposition 3.3]{Tapiero2} for more details.

\begin{figure}[h!]
\begin{center}
\includegraphics[width=10cm]{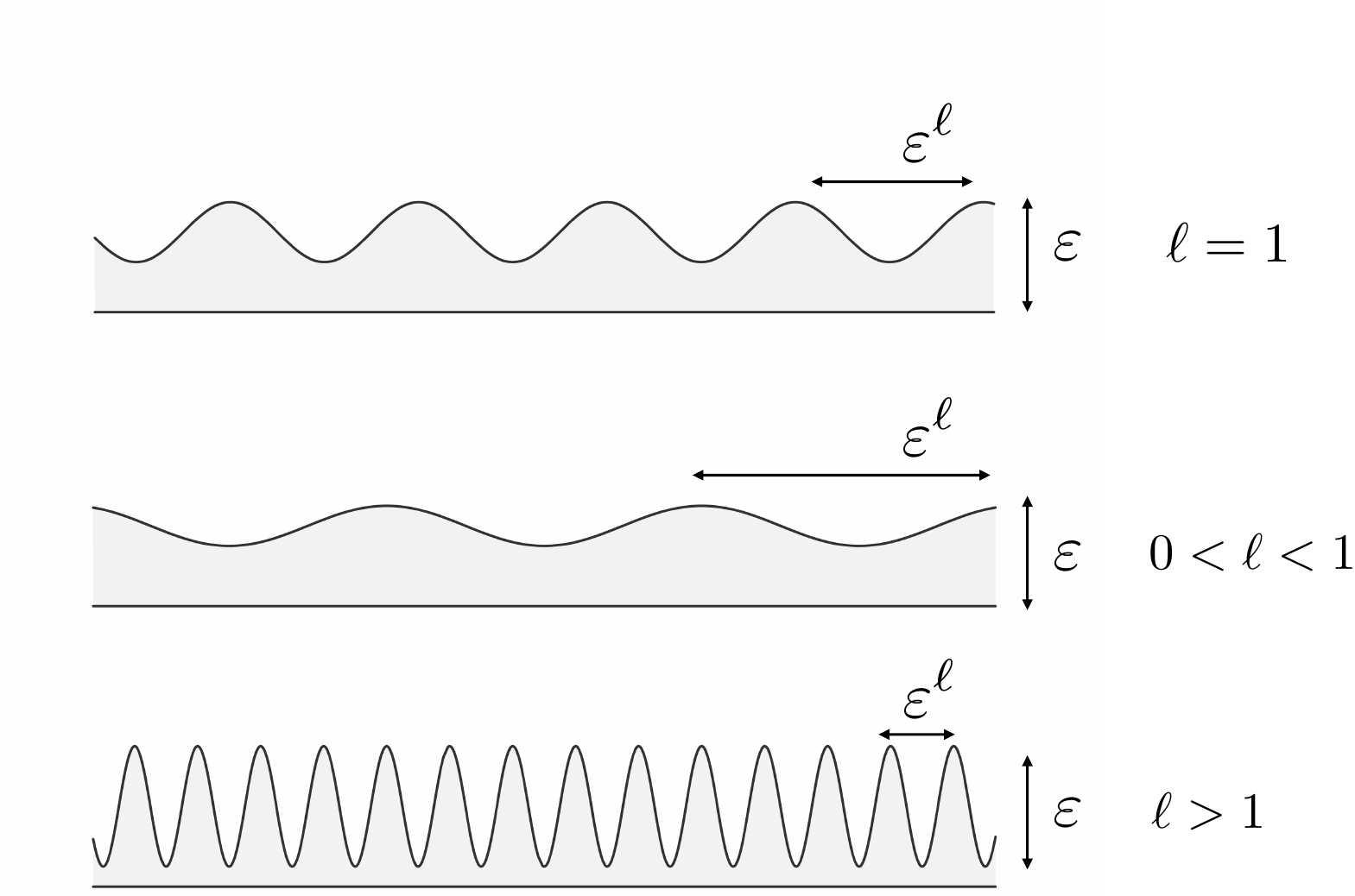}
\end{center}
\caption{View of  {\it Stokes} (top), {\it Reynolds} (middle) and {\it High-frequency} (bottom) roughness regimes.}
\label{fig:omep}
\end{figure}

 On the other hand, engineering practice also stresses the interest of studying the effects of domain irregularities on a thin film flow, see for instance   Dowson \& Higginson \cite{Dowson},   Patir \& Cheng \cite{Patir}, Sullivan \& Sharma \cite{Sullivan} or Tanner \cite{Tanner}.  Thus, the goal becomes in identifying in which way the irregular boundary affects the flow. In this sense, the oscillating boundary is described by two parameters, $\varepsilon$ and $\varepsilon^\ell$, with $\ell\in(0,+\infty)$. The parameter $\varepsilon$ is the thickness of the domain (the distance between the surfaces) and $\ep^\ell$ is the characteristic wavelength of the periodic roughness, see Figure \ref{fig:omep}. This, the function that describes the oscillating boundary is given by $h_\ep(x')=\ep h(x'/\ep^\ell)$, where function $h$ is periodic (see Section \ref{sec:definition} for more detail).

  Employing asymptotic techniques with respect to the domain's thickness, it is showed in Bayada \& Chambat \cite{Bayada_Chambat_1988, Bayada_Chambat}  that depending on the relation between $\ep$ and $\ep^\ell$  there exists three different regimes (see see Figure \ref{fig:omep}):
 \begin{itemize}
 \item[--] {\it Stokes roughness} ($\ell=1$): the roughness oscillation wavelength is of the same order as the thin domain thickness. This models lubrication films where the roughness is quite pronounced relative to the film thickness, such as in bearings or seals with moderate surface textures. The oscillations significantly affect local velocity and pressure distribution inside the lubricant film.
 
 \item[--] {\it Reynolds roughness regime} ($0<\ell<1$): the roughness oscillations have a longer wavelength relative to thickness but still contribute to effective flow properties.  Typical applications are for lubrication layers in devices with smoother surfaces or when the roughness patterns are more spread out, e.g., polymer solutions in thin-film coatings or fluid films in microfluidics with mildly rough channels.
 
 \item[--] {\it High-frequency roughness} ( $\ell>1$):the roughness oscillates on a much smaller scale compared to film thickness, i.e., very rapid boundary oscillations.  Appropriate in highly textured or micro/nano-structured surfaces where surface features are very fine relative to the lubricant film thickness, such as in advanced tribological coatings or engineered superhydrophobic surfaces.  
 \end{itemize} 

 These characteristic regimes are of great interest, for instance they have been also studied in  Benhaboucha {\it et al.} \cite{Benhaboucha}, Fabricius {\it et al.} \cite{Fabricius} and Mikeli\'c \cite{Mikelic2} for Newtonian fluids, in Anguiano \& Su\'arez-Grau  \cite{Anguiano_SG}  and Nakasato \& Pa\v zanin \cite{Nakasato} for power law fluids, Cardone {\it et al.} \cite{Cardone} for Bingham fluids, and  in Pa\v zanin \& Su\'arez-Grau \cite{Pazanin2} and Su\'arez-Grau \cite{grau2} for micropolar fluids.

  The {\it high-frequency roughness} regime ($\ell>1$), due to the highly oscillating boundary, leads to the conclusion that
the velocity is zero in the roughness zone (see e.g. \cite{Anguiano_SG}) so, in this paper, we study the
asymptotic behavior of the solution in the {\it Stokes roughness} ($\ell=1$) and {\it Reynolds roughness} ($0<\ell<1$) regimes.

This paper addresses this gap by rigorously deriving a nonlinear two-dimensional Reynolds model for a quasi-Newtonian fluid governed by the Carreau law flowing through a thin domain with a rapidly oscillating boundary,  which as far as the authors know, it has not yet been studied. Thus, in the spirit of the previous studies,   the following modified Reynolds law (see Section \ref{sec:limitmodel} for more detail) is derived
$${\rm div}_{x'}\left(\mathcal{U}(\nabla_{x'}\widetilde p(x')-{\bf f}'(x'))\right)=0,$$
where the permeability function $\mathcal{U}:\mathbb{R}^2\to \mathbb{R}^2$ has different expressions depending on the {\it Reynolds roughness regime} $0<\ell<1$ (see Section \ref{sec:Reynolds}) or the {\it Stokes roughness regime} $\ell=1$ (see Section \ref{sec:Stokes}).  The novelty lies in the simultaneous treatment of the nonlinear viscosity dependence given by the Carreau model, which adds complexity due to the nonlinear relation between the shear rate and the viscosity; the presence of a slightly rough boundary with small-scale oscillations, whose interplay with the thin-film geometry introduces additional multi-scale challenges in the asymptotic analysis; the rigorous mathematical justification of the homogenization limit as the thickness parameter tends to zero, while capturing effective contributions from the oscillating boundary, and the incorporation of monotonicity arguments and sharp a priori estimates to establish compactness and convergence of the nonlinear terms.

The transition from classical Newtonian to quasi-Newtonian fluids governed by the Carreau law introduces several substantial technical differences. Primarily, the viscosity’s nonlinear dependence on the shear rate replaces the constant viscosity, converting the Stokes system into a nonlinear PDE with a monotone but nonlinear operator. Establishing existence, uniqueness, and compactness for the solutions requires monotonicity arguments and sharp a priori estimates rather than the classical linear elliptic theory. Moreover, the effective two-dimensional Reynolds equation derived in the thin-film asymptotic limit becomes nonlinear and involves integral expressions accounting for the non-Newtonian rheology and the interaction with the rough boundary. The Carreau law’s introduction of zero and infinite shear rate viscosity plateaus regularizes the model but complicates functional dependence. Thus, the nonlinear coupling necessitates more advanced homogenization and unfolding techniques beyond classical approaches utilized in Newtonian flow analyses.

The mathematical analysis is based on deriving the sharp {\it a priori} estimates using a decomposition of the pressure in two pressures, one more regular and depending only on the horizontal variables, and other pressure less regular.  Thanks to the application of an adaptation of the partial unfolding method (only applied to the horizontal variables), and proving the compactness results of the rescaled (in the vertical variable) functions together with monotonicity arguments, we pass to the limit by obtaining reduction of dimension.

By overcoming these difficulties, the derived nonlinear Reynolds equation captures the essential influence of the roughness on the hydrodynamic behavior of quasi-Newtonian fluids in lubrication regimes. This contributes both theoretically and practically, providing a model that can be used in applications involving complex fluids lubricated over rough surfaces, where both nonlinearity of viscosity and surface microstructure play key roles. As stated previously, we have developed the cases of {\it Reynolds }and {\it Stokes roughness regime}. Moreover, we also give some comments concerning the corresponding result for the {\it high-frequency roughness regime} in Section \ref{sec:high}.

 The effective nonlinear two-dimensional Reynolds model obtained in this study not only advances the theoretical understanding of quasi-Newtonian fluid flow under the Carreau law in thin domains with rough boundaries but also has significant practical importance in various engineering and industrial settings. For instance, the results are directly applicable to lubrication problems involving complex fluids such as polymer melts, paints, biological fluids, and other shear-thinning or shear-thickening substances where the classical Newtonian assumption fails. The inclusion of rough boundaries enhances the model's realism, addressing essential micro-scale surface roughness effects frequently encountered in tribology and microfluidic devices. Moreover, this improved understanding can aid in the design and optimization of lubrication systems for biomedical implants, polymer processing equipment, and coatings technology where control over flow behavior near rough interfaces is crucial. The model can also support predictive analysis and simulation in emerging fields such as soft robotics and flexible electronics, where thin non-Newtonian fluid layers intervene in mechanical contacts with irregular surfaces.

 Although the present work focuses on steady-state quasi-Newtonian Stokes flow governed by the Carreau law, the nonsteady or transient regimes are of significant practical and theoretical interest. Time-dependent quasi-Newtonian flows, characterized by evolving shear-dependent viscosities and possibly viscoelastic effects, are critically important in many industrial processes and biological systems. The mathematical analysis of such time-dependent flows generally requires advanced functional frameworks to handle nonlinear and time-dependent viscosity, and several studies have addressed numerical and analytical aspects in this context (see e.g. Baranger \& Najib \cite{Baranger}, Mikeli\'c \cite{Mikelic2, MikelicIntro1}, and recent works by Lu \& Quian \cite{Lu} and Nakasato \& Pa\v zanin  \cite{Nakasato}). Extending the asymptotic analysis and derivation of effective two-dimensional Reynolds-type models to nonsteady quasi-Newtonian fluids remains a challenging and relevant research direction.

Finally, we comment the structure of the paper. In Section \ref{sec:definition} we define the domain and introduce some notation which will be useful in the rest of the paper. In Section \ref{sec:model}, we present the model that will be  analyzed. A priori estimates of velocity and pressure (using a decomposition result for the pressure) are derived in Section \ref{sec:estimates} depending on the type of fluid (pseudoplastic $1<r<2$ or dilatant $r>2$).  A version of the unfolding method to capture the oscillations of the domain is introduced in Section \ref{sec:unfolding}, which gives different estimates depending on the regime $0<\ell<1$ or $\ell>1$ and, in each regime, depending on the type of fluid.  Convergences of velocity and pressure are presented in Section \ref{sec:convergences} and the limit models (main results) are derived in Section \ref{sec:limitmodel}. We finish the paper with a list of references.

\section{Definition of the domain and some notation}\label{sec:definition}   We define the thin domain with a   rough boundary  by
\begin{equation}\label{Omegaep}
\Omega_\varepsilon=\{x=(x',x_3)\in\mathbb{R}^2\times \mathbb{R}\,:\, x'\in \omega,\ 0<x_3< h_\varepsilon(x')\},
\end{equation}
where  $\omega$ a open connected subset of $\mathbb{R}^2$ with Lipschitz boundary and the function $h_\ep(x')= \varepsilon h\left(x'/\varepsilon^\ell\right)$ represents the real distance between the two surfaces. The small parameter $0<\varepsilon\ll 1$ is related to the film thickness and  $0<\varepsilon^\ell\ll 1$ is the wavelength of the roughness. We consider that $\varepsilon^\ell$ is of order greater or equal than $\varepsilon$, i.e. $0<\ell\leq 1$. Thus we consider two cases: the {\it Stokes roughness regime} with $\ell=1$ and the {\it Reynolds roughness regime} $0<\ell<1$ which implies
\begin{equation}\label{RelationAlpha}
\lim_{\ep\to 0}\ep^{1-\ell}=0.
\end{equation}

\noindent Function $h$ is a strictly positive, bounded  and  Lipschitz function defined for $z'$ in $\mathbb{R}^2$, $Z'$-periodic with $Z'=(-1/2,1/2)^2$ the cell of periodicity in $\mathbb{R}^2$, and there exist $h_{\rm min}$ and $h_{\rm max}$ such that
$$0<h_{\rm min}=\min_{z'\in Z'} h(z'),\quad  h_{\rm max}=\max_{z'\in Z'}h(z')\,.$$

\noindent
We define the boundaries of $\Omega_\varepsilon$ as follows
$$\begin{array}{c}
\displaystyle \Gamma_0 =\omega\times\{0\},\quad \Gamma_1^\varepsilon=\left\{(x',x_3)\in\mathbb{R}^2\times\mathbb{R}\,:\, x'\in \omega,\ x_3=  h_\varepsilon(x')\right\},\quad
\displaystyle\Gamma_{\rm lat}^\varepsilon=\partial\Omega_\ep\setminus (\Gamma_0\cup\Gamma_1^\ep).
\end{array}$$

To define the microstructure of the periodicity of the boundary,  we assume that the domain $\omega$ is divided by a mesh of size $\ep^\ell$: for $k'\in\mathbb{Z}^2$, each cell $Z'_{k',\ep^\ell}=\ep^\ell k'+\ep^\ell Z'$. We define $\mathcal{T}_\ep=\{k'\in\mathbb{Z}\,:\, Z'_{k',\ep^\ell}\cap\omega\neq \emptyset\}$. In this setting, there exists an exact finite number of periodic sets $Z'_{k',\ep^\ell}$ such that $k'\in T_\ep$.   Also, we define $Z_{k',\ep^\ell}=Z'_{k',\ep^\ell}\times (0,h(z'))$ and  $Z=Z'\times (0,h(z'))$, which is the reference cell in $\mathbb{R}^3$. We define the boundaries $\widehat \Gamma_0=Z'\times \{0\}$, $\widehat \Gamma_1=Z'\times\{h(z')\}$.

\noindent Applying a dilatation in the vertical variable, i.e. $z_3=x_3/\ep$, we define the following rescaled sets 
 \begin{equation}\label{domains_tilde}\begin{array}{c}
 \displaystyle \widetilde \Omega_\varepsilon=\{(x',z_3)\in\mathbb{R}^2\times \mathbb{R}\,:\, x'\in \omega,\ 0<z_3< h(x'/\ep^\ell)\},\\
 \noame
  \widetilde\Gamma_1^\varepsilon=\{(x',z_3)\in\mathbb{R}^2\times \mathbb{R}\,:\, x'\in \omega,\ z_3=h(x'/\ep^\ell)\}\quad \widetilde \Gamma_{\rm lat}^\varepsilon=\partial\widetilde \Omega_\ep\setminus (\widetilde \Gamma_0\cup\widetilde \Gamma_1^\ep).
  \end{array}
  \end{equation}
The quantity  $h_{\rm max}$ allows us to define:
\begin{itemize}
\item[--] The extended sets $Q_\ep=\omega\times (0,\ep h_{\rm max})$,  $\Omega=\omega\times  (0, h_{\rm max})$  and $\Gamma_1=\omega\times \{h_{\rm max}\}$.
\item[--]  The extended cube  $\widetilde Q_{k',\varepsilon^\ell}=Z'_{k',\varepsilon^\ell}\times (0, h_{\rm max})$ for $k'\in\mathbb{Z}^2$.

\end{itemize}

\noindent We denote by $C$ a generic constant which can change from line to line. Moreover, $O_\ep$ denotes a  generic quantity, which can change from line to line, devoted to tend to zero when $\ep\to 0$.\\

\noindent In the sequel, we introduce the following notation. Let us consider a vectorial function ${\bf v} =( {\bf v} ', v_{3})$ with ${\bf v}'=(v_1, v_2)$ and $\varphi$ a scalar function, both defined in $\Omega_\varepsilon$. Then, introduce the operator $\Delta$ and $\nabla$ by
$$\begin{array}{c}
 \displaystyle \Delta{\bf v}=\Delta_{x'}{\bf v} +\partial_{x_3}^2 {\bf v},\quad {\rm div}({\bf v} )={\rm div}_{x'}({\bf v}')+\partial_{x_3}v_3,\quad \nabla \varphi=(\nabla_{x'} \varphi,  \partial_{x_3}\varphi)^t,\\
 \end{array}$$
and
$\mathbb{D}:\mathbb{R}^3\to \mathbb{R}^3_{\rm sym}$   the symmetric part of the velocity gradient, that is
$$\mathbb{D}[{\bf v}]={1\over 2}(D{\bf v}+(D{\bf v})^t)=\left(\begin{array}{ccc}
\partial_{x_1}v_1 &   {1\over 2}(\partial_{x_1}v_2 + \partial_{x_2}v_1) &   {1\over 2}(\partial_{x_3}v_1 + \partial_{x_1}v_3)\\
\noame
 {1\over 2}(\partial_{x_1}v_2 + \partial_{x_2}v_1) & \partial_{x_2}v_2 &   {1\over 2}(\partial_{x_3}v_2 + \partial_{x_2}v_3)\\
 \noame
 {1\over 2}(\partial_{x_3}v_1 + \partial_{x_1}v_3)&   {1\over 2}(\partial_{x_3}v_2 + \partial_{x_2}v_3)& \partial_{x_3}v_3
\end{array}\right).$$

\noindent Moreover, for $\widetilde {\bf  v} =(\widetilde {\bf  v}', \widetilde  v_{3})$  a vector function and $\widetilde \varphi$ a scalar function, both defined in $\widetilde\Omega_\ep$, obtained from ${\bf v}$ and $\varphi$ after a dilatation in the vertical variable, respectively, we will use the following operators
 $$\begin{array}{c}
 \displaystyle \Delta_{\varepsilon}  \widetilde  {\bf v} =\Delta_{x'} \widetilde  {\bf v} + \varepsilon^{-2}\partial_{z_3}^2 \widetilde  {\bf v} ,\quad   \displaystyle\Delta_{ \varepsilon}\widetilde \varphi=\Delta_{x'}\widetilde \varphi+ \ep^{-2}\partial^2_{z_3}\widetilde \varphi,\\
  \noame
  (D_{\ep}\widetilde {\bf v})_{ij}=\partial_{x_j}\widetilde {v}_i\ \hbox{ for }\ i=1,2,3,\ j=1,2,\quad   (D_{\ep}\widetilde {\bf v})_{i3}=\ep^{-1}\partial_{z_3}\widetilde {v}_i \ \hbox{ for }\ i=1,2,3,\\
  \noame
  \nabla_{\ep}\widetilde\varphi=(\nabla_{x'}\widetilde \varphi, \ep^{-1}\partial_{z_3}\widetilde\varphi)^t,\quad {\rm div}_{\varepsilon}(\widetilde  {\bf   v})={\rm div}_{x'}(\widetilde {\bf v}')+\ep^{-1}{\partial_{z_3}}\widetilde  v_{3},
  \end{array}$$
 Moreover, we define $\mathbb{D}_{\ep}[\widetilde  {\bf v}]$ as follows
 $$\mathbb{D}_{\ep}[\widetilde  {\bf v}]=\mathbb{D}_{x'}[\widetilde  {\bf v}]+\ep^{-1}\partial_{z_3}[\widetilde  {\bf v}]=\left(\begin{array}{ccc}
\partial_{x_1}\widetilde  v_1 &   {1\over 2}(\partial_{x_1}\widetilde  v_2 + \partial_{x_2}\widetilde  v_1) &   {1\over 2} (\partial_{x_1}\widetilde  v_3+_\ep^{-1}\partial_{z_3}\widetilde  v_1)\\
\noame
 {1\over 2}(\partial_{x_1}\widetilde  v_2 + \partial_{x_2}\widetilde  v_1) & \partial_{x_2}\widetilde  v_2 &   {1\over 2} (\partial_{x_2}\widetilde  v_3+_\ep^{-1}\partial_{z_3}\widetilde  v_2)\\
 \noame
 {1\over 2} (\partial_{x_1}\widetilde  v_3+\ep^{-1}\partial_{z_3}\widetilde  v_1)&   {1\over 2} (\partial_{x_2}\widetilde  v_3+\ep^{-1}\partial_{z_3}\widetilde  v_2)& \ep^{-1}\partial_{z_3}\widetilde  v_3
\end{array}\right),$$
 where $\mathbb{D}_{x'}[\widetilde  {\bf v}]$ and $\partial_{z_3}[\widetilde  {\bf v}]$ are defined by
\begin{equation}\label{def_der_sym_1}
  \mathbb{D}_{x'}[{\bf v}]=\left(\begin{array}{ccc}
\partial_{x_1}v_1 &   {1\over 2}(\partial_{x_1}v_2 + \partial_{x_2}v_1) &   {1\over 2} \partial_{x_1}v_3\\
\noame
 {1\over 2}(\partial_{x_1}v_2 + \partial_{x_2}v_1) & \partial_{x_2}v_2 &   {1\over 2} \partial_{x_2}v_3\\
 \noame
 {1\over 2} \partial_{x_1}v_3&   {1\over 2} \partial_{x_2}v_3& 0
\end{array}\right),
 \  \partial_{z_3}[{\bf v}]=\left(\begin{array}{ccc}
0 &   0&   {1\over 2} \partial_{z_3}v_1\\
\noame
 0& 0 &   {1\over 2} \partial_{z_3}v_2\\
 \noame
 {1\over 2} \partial_{z_3}v_1&   {1\over 2} \partial_{z_3}v_2& \partial_{z_3}v_3
\end{array}\right).
\end{equation}

\noindent We also define the following operators applied to ${\bf v}'$:
\begin{equation}\label{def_der_sym_2}
  \mathbb{D}_{x'}[{\bf v}']=\left(\begin{array}{ccc}
\partial_{x_1}v_1 &   {1\over 2}(\partial_{x_1}v_2 + \partial_{x_2}v_1) &   0\\
\noame
 {1\over 2}(\partial_{x_1}v_2 + \partial_{x_2}v_1) & \partial_{x_2}v_2 &  0\\
 \noame
0&   0& 0
\end{array}\right),
 \quad \partial_{z_3}[{\bf v}']=\left(\begin{array}{ccc}
0 &   0&   {1\over 2} \partial_{z_3}v_1\\
\noame
 0& 0 &   {1\over 2} \partial_{z_3}v_2\\
 \noame
 {1\over 2} \partial_{z_3}v_1&   {1\over 2} \partial_{z_3}v_2& 0
\end{array}\right).
\end{equation}

\section{Model problem}\label{sec:model} With homogeneous Dirichlet boundary conditions, the flow of velocity ${\bf u}_\ep=({\bf u}'_\ep(x), u_{3,\ep}(x))$ and pressure $p_\ep=p_\ep(x),$ at a point $x\in \Omega_\ep$, is suppose to be ruled by the following  Stokes system
\begin{equation}\label{system_1}
\left\{\begin{array}{rl}
\displaystyle -\varepsilon{\rm div}(\eta_r(\mathbb{D}[{\bf u}_\ep])\mathbb{D}[{\bf u}_\ep])+\nabla p_\ep={\bf f} & \hbox{in}\ \Omega_\varepsilon,\\
\noame
{\rm div}( {\bf u}_\varepsilon)=0& \hbox{in}\ \Omega_\varepsilon\\
\noame
\displaystyle {\bf u}_\varepsilon=0 & \hbox{on }\partial\Omega_\varepsilon.
\end{array}\right.
\end{equation}
Here, $\eta_r$ is the Carreau law defined  in (\ref{Carreaulaw}).  In the momentum equation (\ref{system_1}), according to \cite[Case $\gamma=1$]{Tapiero2}, the viscosity $\eta_r$ is affected by a factor $\ep$, which can  be interpreted as a scaling of the low-shear-rate and high-shear-rate limit viscosities $\eta_0, \eta_\infty$ appearing in the Carreau law (\ref{Carreaulaw}) by $\ep$.  The source term ${\bf f}$ is of the form
\begin{equation}\label{fassump}
{\bf f} (x)=({\bf f}'(x'),0)\quad \hbox{with}\quad {\bf f}'\in L^\infty(\omega)^2.
\end{equation}
Notice that the assumptions of neglecting the vertical component  of the exterior force and  the independence of the vertical variable  are usual when dealing with fluids through thin domains (see \cite{Tapiero2} for more details).
\\
\begin{remark} Under previous assumptions, for every fixed positive $\varepsilon$, the classical theory (see for instance \cite{Tapiero2,Bourgeat1,Lions2}), gives the existence of a unique weak  solution $({\bf u}_\varepsilon, p_\varepsilon)\in H^1_0(\Omega_\varepsilon)^3\times L^{2}_0(\Omega_\varepsilon)$, for $1<r< 2$, and $({\bf u}_\varepsilon, p_\varepsilon)\in W^{1,r}_0(\Omega_\varepsilon)^3\times L^{r'}_0(\Omega_\varepsilon)$ with $1/r+1/r'=1$, for $r> 2$, where  $L^{2}_0$ (respectively $L^{r'}_0$) is the space of functions of $L^2$ (respectively $L^{r'}$) with zero mean value.
\end{remark}
To study the asymptotic behavior of the solutions ${\bf u}_\ep$ and $p_\ep$ when $\ep$  tends to zero,  we use the rescaling  
\begin{equation}\label{dilatacion}
z_3={x_3\over \ep}\,,
\end{equation}
 to have the functions defined in  $\widetilde\Omega_\ep$, which is defined in (\ref{domains_tilde}). Using the change of variables (\ref{dilatacion}) in  the model problem, we get the rescaled  Stokes system
\begin{equation}\label{system_1_dil}
\left\{\begin{array}{rl}
\displaystyle -\varepsilon {\rm div}_\varepsilon(\eta_r(\mathbb{D}_\varepsilon[ \widetilde {\bf u}_\ep])\mathbb{D}_\varepsilon[\widetilde {\bf u}_\ep])+\nabla_\varepsilon \widetilde p_\ep={\bf f} & \hbox{in}\ \widetilde \Omega_\varepsilon,\\
\noame
{\rm div}_\varepsilon( \widetilde {\bf u}_\varepsilon)=0& \hbox{in}\ \widetilde \Omega_\varepsilon\\
\noame
\displaystyle \widetilde {\bf u}_\varepsilon=0 & \hbox{on }\partial\widetilde \Omega_\varepsilon,
\end{array}\right.
\end{equation}

\noindent where the unknown functions in the above system are given by ${\bf \widetilde u}_\varepsilon(x',z_3)={\bf u}_\varepsilon(x',\varepsilon z_3)$, $ \widetilde p_\varepsilon(x',z_3)=p_\varepsilon(x',\varepsilon z_3)$  for a.e. $(x',z_3)\in \widetilde\Omega_\varepsilon$, and the operators ${\rm div}_\ep$, $\mathbb{D}_\ep$ and $\nabla_\ep$ are defined in Section \ref{sec:definition}.

Next, we derive the sharp {\it a priori} estimates, which will allow us to later prove  compactness results of the rescaled functions.

\section{{\it A priori} estimates} \label{sec:estimates}

 To derive the desired estimates for velocity,  let us recall and prove some well-known technical results (see for instance \cite{Anguiano_SG}).
\begin{lemma}[Poincar\'e's  and Korn's inequalities]\label{Poincare_lemma} For every $\varphi\in W^{1,q}_0(\Omega_\varepsilon)^3$, $1\leq q<+\infty$, it holds
\begin{equation}\label{Poincare}
\|\varphi\|_{L^q(\Omega_\varepsilon)^3}\leq C\varepsilon\|D \varphi\|_{L^q(\Omega_\varepsilon)^{3\times 3}},
\quad
 \|D\varphi\|_{L^q(\Omega_\varepsilon)^{3\times 3}}\leq C\|\mathbb{D}[\varphi]\|_{L^q(\Omega_\varepsilon)^{3\times 3}}.
\end{equation}
Moreover, from the change of variables (\ref{dilatacion}),    for every $\widetilde\varphi\in W^{1,q}_0(\widetilde\Omega_\varepsilon)^3$, it holds  the following rescaled estimate
\begin{equation}\label{Poincare2}
\|\widetilde \varphi\|_{L^q(\widetilde \Omega_\varepsilon)^3}\leq C\varepsilon\|D_{\ep} \widetilde \varphi\|_{L^q(\widetilde \Omega_\varepsilon)^{3\times 3}},\quad
 \|D_\varepsilon\widetilde \varphi\|_{L^q(\widetilde \Omega\varepsilon)^{3\times 3}}\leq C\|\mathbb{D}_\varepsilon[\widetilde \varphi]\|_{L^q(\widetilde \Omega_\varepsilon)^{3\times 3}}.
\end{equation}
\end{lemma}

\begin{lemma}\label{lemma_estimates} Let ${\bf u}_\ep$ be the solution of (\ref{system_1}) and $\widetilde {\bf u}_\ep$  the solution of  (\ref{system_1_dil}).  We have the following estimates depending on the value of  $r$:
\begin{itemize}
\item  If $1<r<+\infty$, $r\neq 2$, the following estimates hold
\begin{equation}\label{estim_sol1}
\displaystyle\|{\bf   u}_\varepsilon\|_{L^2(\Omega_\varepsilon)^3}\leq C\varepsilon^{3\over 2}, \quad\displaystyle
\|D{\bf  u}_\varepsilon\|_{L^2(\Omega_\varepsilon)^{3\times 3}}\leq C\ep^{1\over 2},\quad\displaystyle
\|\mathbb{D}[{\bf u}_\varepsilon]\|_{L^2(\Omega_\varepsilon)^{3\times 3}}\leq C\ep^{1\over 2},
\end{equation}
\begin{equation}\label{estim_sol_dil1}
\|{\bf \widetilde u}_\varepsilon\|_{L^2(\widetilde \Omega_\varepsilon)^3}\leq C\varepsilon, \quad\displaystyle
\|D_{\varepsilon} {\bf  \widetilde u}_\varepsilon\|_{L^2(\widetilde\Omega_\varepsilon)^{3\times 3}}\leq C,\quad\displaystyle
\|\mathbb{D}_{\varepsilon} [{\bf  \widetilde u}_\varepsilon]\|_{L^2(\widetilde\Omega_\varepsilon)^{3\times 3}}\leq C.
\end{equation}
 \item  If $r>2$ the following estimates hold
\begin{equation}\label{estim_sol1_r}
\displaystyle
\|{\bf u}_\varepsilon\|_{L^r( \Omega_\varepsilon)^3}\leq C\varepsilon^{1+{1\over r}}, \quad\displaystyle
\|D {\bf  u}_\varepsilon\|_{L^r(\Omega_\varepsilon)^{3\times 3}}\leq C\ep^{1\over r},\quad\displaystyle
\|\mathbb{D} [{\bf  u}_\varepsilon]\|_{L^r(\Omega_\varepsilon)^{3\times 3}}\leq C\ep^{1\over r},
\end{equation}
\begin{equation}\label{estim_sol_dil1_r}
\|{\bf \widetilde u}_\varepsilon\|_{L^r(\widetilde \Omega_\varepsilon)^3}\leq C\varepsilon, \quad\displaystyle
\|D_{\varepsilon} {\bf  \widetilde u}_\varepsilon\|_{L^r(\widetilde\Omega_\varepsilon)^{3\times 3}}\leq C,\quad\displaystyle
\|\mathbb{D}_{\varepsilon} [{\bf  \widetilde u}_\varepsilon]\|_{L^r(\widetilde\Omega_\varepsilon)^{3\times 3}}\leq C.
\end{equation}
\end{itemize}

\end{lemma}
\begin{proof}
Multiplying (\ref{system_1}) by $ {\bf u}_\varepsilon$, integrating over $\Omega_\varepsilon$ and taking into account  that ${\rm div}({\bf u}_\varepsilon)=0$ in $\Omega_\varepsilon$,  we get
\begin{equation}\label{form_var_estim}\varepsilon (\eta_0-\eta_\infty)\int_{\Omega_\varepsilon}\left(1+\lambda|\mathbb{D}[\ {\bf u}_\varepsilon]|^2\right)^{{r\over 2}-1}|\mathbb{D}[  {\bf u}_\varepsilon]|^2dx+\varepsilon  \eta_\infty\int_{ \Omega_\varepsilon}|\mathbb{D}[ {\bf u}_\varepsilon]|^2dx= \int_{\Omega_\varepsilon}{\bf f}'\cdot \widetilde {\bf u}_\varepsilon'\,dx.
\end{equation}
We divide the proof in two steps. First, we derive estimates (\ref{estim_sol1}) for every $r>1$ and then, for $r>2$, we establish estimates (\ref{estim_sol1_r}). Estimates (\ref{estim_sol_dil1}) and (\ref{estim_sol_dil1_r}) are consequences of the change of variables (\ref{dilatacion}) to estimates (\ref{estim_sol1})  and (\ref{estim_sol1_r}), respectively.\\

{\it Step 1}. We consider $r>1$, $r\neq 2$. Taking into account that $\eta_0>\eta_{\infty}$, and $\lambda>0$, we have
$$\varepsilon(\eta_0-\eta_\infty)\int_{\Omega_\varepsilon}\left(1+\lambda|\mathbb{D}[{\bf u}_\varepsilon]|^2\right)^{{r\over 2}-1}|\mathbb{D}[ {\bf u}_\varepsilon]|^2dx\geq 0.
$$
From  the Cauchy-Schwarz inequality and the assumption on ${\bf f}'$ given in (\ref{fassump}), we deduce from (\ref{form_var_estim}) that
$$\varepsilon \eta_\infty\| \mathbb{D}[ {\bf u}_\varepsilon]\|_{L^2(\Omega_\varepsilon)^{3\times 3}}^2\leq C\ep^{1\over 2} \|  {\bf u}_\varepsilon\|_{L^2( \Omega_\varepsilon)^3}.$$
Applying the Poincar\'e  and Korn inequalities (\ref{Poincare}) to the right-hand side, we get (\ref{estim_sol1})$_3$. Finally, applying once again (\ref{Poincare})  yields (\ref{estim_sol1})$_1$ and (\ref{estim_sol1})$_2$.\\

{\it Step 2}. Assume that $r>2$. The idea is now to estimate the first integral in~\eqref{form_var_estim}. Taking into account that 
$$\varepsilon  \eta_\infty\int_{\Omega_\varepsilon}|\mathbb{D}[{\bf u}_\varepsilon]|^2dx\geq 0,$$
 hence, \eqref{form_var_estim} and the Cauchy-Schwarz inequality imply
\begin{equation}\label{form_var_estim2}\varepsilon  (\eta_0-\eta_\infty)\int_{\Omega_\varepsilon}\left(1+\lambda|\mathbb{D}[ {\bf u}_\varepsilon]|^2\right)^{{r\over 2}-1}|\mathbb{D}[ {\bf u}_\varepsilon]|^2dx\leq  C\ep^{1\over 2}\| {\bf u}_\varepsilon\|_{L^2(\Omega_\varepsilon)^3}.
\end{equation}
Noticing that
$$\varepsilon  \lambda^{r-2\over 2} (\eta_0-\eta_\infty)\int_{\Omega_\varepsilon}|\mathbb{D}[ {\bf u}_\varepsilon]|^rdx\leq \varepsilon  (\eta_0-\eta_\infty)\int_{\Omega_\varepsilon}\left(1+\lambda|\mathbb{D}[ {\bf u}_\varepsilon]|^2\right)^{{r\over 2}-1}|\mathbb{D}[{\bf u}_\varepsilon]|^2dx,$$
and applying the Poincar\'e and Korn inequalities (\ref{Poincare}) to the right-hand side of (\ref{form_var_estim2}), we get
\begin{equation}\label{estim_du}\ep\|\mathbb{D}[ {\bf u}_\varepsilon]\|^r_{L^r( \Omega_\varepsilon)^{3\times 3}}\leq C \ep^{1+{1\over 2}}\|\mathbb{D}[{\bf u}_\varepsilon]\|_{L^2(\Omega_\varepsilon)^{3\times 3}}.
\end{equation}
On the one hand, applying estimate (\ref{estim_sol1})$_3$, we deduce
$$\|\mathbb{D}[ {\bf u}_\varepsilon]\|_{L^r(\Omega_\varepsilon)^{3\times 3}}\leq C\ep^{1\over r}.$$
On the other hand, similarly to (\ref{form_var_estim2}), we also deduce
\begin{equation}\label{form_var_estim21212}\varepsilon  (\eta_0-\eta_\infty)\int_{\Omega_\varepsilon}\left(1+\lambda|\mathbb{D}[ {\bf u}_\varepsilon]|^2\right)^{{r\over 2}-1}|\mathbb{D}[ {\bf u}_\varepsilon]|^2dx\leq  C\ep^{1\over r'}\| {\bf u}_\varepsilon\|_{L^r(\Omega_\varepsilon)^3}.
\end{equation}
Proceeding similarly, from the continuity of the embedding $L^r(\Omega_\varepsilon)\hookrightarrow L^{2}(\Omega_\varepsilon)$ in (\ref{estim_du}), we also have
$$\ep\|\mathbb{D}[ {\bf u}_\varepsilon]\|^r_{L^r(\Omega_\varepsilon)^{3\times 3}}\leq C\ep^{1+{1\over r'}}\|\mathbb{D}[{\bf u}_\varepsilon]\|_{L^r(\Omega_\varepsilon)^{3\times 3}},$$
which also gives $\|\mathbb{D}[ {\bf u}_\varepsilon]\|_{L^r(\Omega_\varepsilon)^{3\times 3}}\leq C\ep^{1\over r}.$
Finally, from the Poincar\'e and Korn inequalities (\ref{Poincare}), we derive the remaining estimates.

\end{proof}

\begin{remark}[Extension of $\widetilde{\bf u}_\ep$ to $\Omega=\omega\times  (0, h_{\rm max})$] From the boundary conditions satisfied by ${\bf \widetilde u}_\varepsilon$ on the top boundary $\widetilde \Gamma_1^\ep$,  we can extend it  by zero in $\Omega\setminus \widetilde{\Omega}^{\varepsilon}$ and we denote the extension by  the same symbol. As consequence, the extended velocity satisfies the same estimates given in Lemma \ref{lemma_estimates}.
\end{remark}

Next, we give a version of a decomposition result which decomposes the pressure $p_\ep$  in two pressures $p_\ep^0$ and $p^1_\ep$, see \cite[Proposition 3.6]{Pazanin2} or \cite[Proposition 2.13]{Nakasato} for more details. This results play an important role to identify the limit problem.
\begin{proposition}\label{prop_decom} Consider $q=\max\{2,r\}$ and $q'$ the conjugate exponent of $q$, such that, $1/q+1/q'=1$. Then, the pressure $p_\ep\in L^{q'}_0(\Omega^\ep)$, solution to (\ref{system_1}), can be decomposed as
\begin{equation}\label{decomposition}
p_\ep=p^\ep_0+p^\ep_1,
\end{equation}
where $p_0^\ep\in W^{1,q'}(\omega)$, which is independent of $x_3$, and $p^\ep_1\in L^{q'}(\Omega_\ep)$. Moreover, the following estimate holds
$$\ep^{1+{1\over q'}}\|p_0^\ep\|_{W^{1,q'}(\omega)}+\|p^\ep_1\|_{L^{q'}(\Omega_\ep)}\leq C\|\nabla p_\ep\|_{W^{-1,q'}(\Omega_\ep)^3},$$
that is
\begin{equation}\label{estim_decom_pressure}
\|p_0^\ep\|_{W^{1,q'}(\omega)}\leq C\ep^{-{1+q'\over q'}}\|\nabla p_\ep\|_{W^{-1,q'}(\Omega_\ep)^3},\quad \|p^\ep_1\|_{L^{q'}(\Omega_\ep)}\leq C\|\nabla p_\ep\|_{W^{-1,q'}(\Omega_\ep)^3}.
\end{equation}

\end{proposition}
Now, we give the estimates for pressures.
\begin{lemma}\label{lemma_est_P} Let $p_\varepsilon$ the pressure solution of (\ref{system_1}) and the decomposition (\ref{decomposition}). Then, we have the following estimates depending on the value of $r$:
\begin{itemize}
\item  if $1<r< 2$, then
\begin{equation}\label{esti_P}
\begin{array}{l}\displaystyle
\|p^\ep_0\|_{H^1(\omega)}\leq C,\quad \|p^\ep_1\|_{L^{2}(\Omega_\ep)}\leq C\ep^{{3\over 2}},\quad \|\widetilde p^\ep_1\|_{L^{2}(\widetilde \Omega_\ep)}\leq C\ep,
\end{array}\end{equation}
\item if $r>2$, then 
\begin{equation}\label{esti_P_r}
\|p^\ep_0\|_{W^{1,r'}(\omega)}\leq C,\quad \|p^\ep_1\|_{L^{r'}(\Omega_\ep)}\leq C\ep^{1+{1\over r'}},\quad \|\widetilde p^\ep_1\|_{L^{r'}(\widetilde \Omega_\ep)}\leq C\ep,
\end{equation}
\end{itemize}
where $\widetilde p_1^\ep$ is obtained from $p^\ep_1$ by using the change of variables (\ref{dilatacion}).
\end{lemma}
\begin{proof} To prove estimates (\ref{esti_P}) and (\ref{esti_P_r}), it is enough to obtain the following estimate
\begin{equation}\label{estim_nabla_p}
\|\nabla p_\ep\|_{W^{-1,q'}(\Omega_\ep)^3}\leq C\ep^{1+{1\over q'}},
\end{equation}
with $q=\max\{2,r\}$ and $q'$ the conjugate exponent of $q$, such that, $1/q+1/q'=1$.

Let us prove the  estimate  \eqref{estim_nabla_p}. For $\varphi\in W^{1,q}_0(\Omega_\ep)$, $q=\max\{2,r\}$ and $q'$ its conjugate exponent, the weak formulation of the Stokes system~\eqref{system_1} is
\begin{equation}\label{Fep2}
\begin{array}{rl}
\displaystyle
\langle\nabla p_\ep,\varphi\rangle_{W^{-1,q'}(\Omega_\ep),W^{1,q}_0(\Omega_\ep)}=&\displaystyle -\varepsilon(\eta_0-\eta_\infty)\int_{\Omega_\varepsilon}(1+\lambda|\mathbb{D}[{ u}_\varepsilon]|^2)^{{r\over 2}-1}\mathbb{D}[ { u}_\varepsilon]:\mathbb{D}[ \varphi]\,dx\\
\noame
&\displaystyle-\varepsilon\eta_\infty\int_{\Omega_\varepsilon}\mathbb{D}[ {\bf u}_\varepsilon]:\mathbb{D}[ \varphi]\,dx+\int_{\Omega_\varepsilon}{ f}'\cdot \varphi'dx.
\end{array}
\end{equation}
First, consider the case $1<r<2$, hence $q=q'=2$. Since  $(1+\lambda|\mathbb{D} [ { u}_\ep |^2)^{{r\over 2}-1}\leq 1$,  applying  Cauchy-Schwarz's inequality and using estimate (\ref{estim_sol1}), we obtain
\begin{equation}\label{Step2_41212}
\begin{array}{rl}
\displaystyle \left|\ep(\eta_0-\eta_\infty)\int_{\Omega_\ep}(1+\lambda|\mathbb{D} [ {\bf u}_\ep]|^2)^{{r\over 2}-1}\mathbb{D} [ {\bf u}_\ep]:\mathbb{D} [    \varphi]\,dx\right| \leq & \displaystyle \ep (\eta_0-\eta_\infty)\int_{\Omega_\ep}|\mathbb{D} [ {\bf u}_\ep]||\mathbb{D} [    \varphi] |dx\\
\noame
\leq &\displaystyle \ep (\eta_0-\eta_\infty)\|\mathbb{D} [  {\bf u}_\ep]\|_{L^2( \Omega_\ep)^{3\times 3}}\|\mathbb{D} [  \varphi]\|_{L^2( \Omega_\ep)^{3\times 3}}\\
\noame
\leq &\displaystyle  C\ep^{3\over 2}\| \varphi\|_{H^1_0(\Omega_\ep)^3}.
\end{array}
\end{equation}
Similarly, we obtain
\begin{equation}\label{Step2_4} 
\displaystyle \left|\ep  \eta_\infty\int_{ \Omega_\ep} \mathbb{D} [ {\bf u}_\ep]:\mathbb{D} [  \varphi]\,dx\right| \leq  C\ep^{3\over 2}\| \varphi\|_{H^1_0(\Omega_\ep)^3}.
\end{equation}
Because ${ f}'={ f}'(x')$ is in $L^\infty(\omega)^2$  and using the Poincar\'e inequality (\ref{Poincare}), we obtain
\begin{equation}\label{Step2_51212} 
\displaystyle \left| \int_{ \Omega_\ep} { f}'\cdot   \varphi'\,dx\right| \leq  C\ep^{1\over 2}\|\varphi\|_{L^2(\Omega_\ep)^3}\leq C\ep^{3\over 2}\| D\varphi\|_{L^2(\Omega_\ep)^{3\times 3}}\leq C\ep^{3\over 2}\| \varphi\|_{H^1_0(\Omega_\ep)^3}.
\end{equation}
Coming back to expression (\ref{Fep2}) with $q'=2$, we deduce from (\ref{Step2_41212})--(\ref{Step2_51212}) the estimate
$$\|\nabla p_\ep\|_{H^{-1}(\Omega_\ep)^3}\leq C\ep^{3\over 2},$$ which is (\ref{estim_nabla_p}) for $q=q'=2$. This finishes the proof in  the case $1<r<2$.\\
 
Now, let us consider the case $r> 2$, which implies $q=r$ and $q'=r'$. Recall that $\Omega$ is the fixed domain. Since $r>2$, $L^r(\Omega)$ is continuously embedded in $L^2(\Omega)$, and since the rescaled function $\widetilde \varphi$ is extended by zero to $\Omega$, there exists a constant $C>0$ such that for every $\ep>0$ and $\varphi\in W^{1,r}_0(\Omega_\ep)$,
 	\begin{align*}
 		\|\mathbb{D}_\ep[ \widetilde \varphi]\|_{L^2(\widetilde \Omega_\ep)^{3\times 3}}& = \|\mathbb{D}_\ep[ \widetilde \varphi]\|_{L^2(\Omega)^{3\times 3}}    \leq C \|\mathbb{D}_\ep[ \widetilde \varphi]\|_{L^r(\Omega)^{3\times 3}} = C \|\mathbb{D}_\ep[ \widetilde \varphi]\|_{L^r(\widetilde \Omega_\ep)^{3\times 3}}.
 	\end{align*}
 	Taking into account that for $s\in \{2,r\}$,
 	\[
 	\|\mathbb{D}[\varphi]\|_{L^s(\Omega_\ep)^{3\times 3}}=\ep^{{1\over s}}\|\mathbb{D}_\ep[ \widetilde \varphi]\|_{L^s(\widetilde \Omega_\ep)^{3\times 3}},
 	\]
 	we deduce that for every $\varphi\in W^{1,r}_0(\Omega_\ep)$,
 	\begin{equation}\label{ScaledContinuousLrL2}
 	\|\mathbb{D}[\varphi]\|_{L^2( \Omega_\ep)^{3\times 3}} \leq C\ep^{\frac{1}{2}-\frac{1}{r}} \|\mathbb{D}[\varphi]\|_{L^r( \Omega_\ep)^{3\times 3}}.
 	\end{equation}
 	Using~\eqref{ScaledContinuousLrL2}, H\"older's inequality, and the inequality $(1+X)^{\alpha}\leq C(1+X^\alpha)$ (which is valid for $X\geq 0$, and $\alpha>0$), we obtain
 	$$\begin{array}{l}
 		\displaystyle
 		\int_{ \Omega_\ep}\left|(1+\lambda|\mathbb{D}[  {\bf u}_\ep]|^2)^{{r\over 2}-1}\mathbb{D}[  {\bf u}_\ep]:\mathbb{D}[   \varphi]\right|dx\\
 		\noame
 		\displaystyle
 		\leq C\left(
 		\int_{ \Omega_\ep}|\mathbb{D} [ {\bf u}_\ep]||\mathbb{D} [  \varphi]|\,dx +\int_{  \Omega_\ep}|\mathbb{D} [ {\bf u}_\ep]|^{r-1}|\mathbb{D} [ \varphi]|\,dx 
 		\right)\\
 		\noame
 		\displaystyle
 		\leq C\left(
 		\|\mathbb{D}[ {\bf u}_\ep]\|_{L^2( \Omega_\ep)^{3\times 3}}\|\mathbb{D} [ \varphi]\|_{L^2(\Omega_\ep)^{3\times 3}}+\|\mathbb{D}[ {\bf u}_\ep]\|^{r-1}_{L^r( \Omega_\ep)^{3\times 3}}\|\mathbb{D}[  \varphi]\|_{L^r(  \Omega_\ep)^{3\times 3}}
 		\right)\\
 		\noame
 		\displaystyle\leq C\left( \ep^{\frac{1}{2}-\frac{1}{r}} 		\|\mathbb{D}[  {\bf u}_\ep]\|_{L^2( \Omega_\ep)^{3\times 3}}+\|\mathbb{D} [{\bf  u}_\ep]\|^{r-1}_{L^r(\Omega_\ep)^{3\times 3}}
 		\right)\|\mathbb{D}[   \varphi]\|_{L^r(\Omega_\ep)^{3\times 3}}.
 	\end{array}$$
 	By~\eqref{Poincare} with $q=r$, $\|\mathbb{D}[   \varphi]\|_{L^r(\Omega_\ep)^{3\times 3}}\leq C\|\varphi\|_{W^{1,r}_0(\Omega_\ep)}$. Combining the previous estimate with the third estimates in~\eqref{estim_sol1} and \eqref{estim_sol1_r}, we obtain
 $$\begin{array}{rl}
	\displaystyle \left|  \ep (\eta_0-\eta_{\infty})	\int_{ \Omega_\ep}\left|(1+\lambda|\mathbb{D}[  {\bf u}_\ep]|^2)^{{r\over 2}-1}\mathbb{D}[  {\bf u}_\ep]:\mathbb{D}[   \varphi]\right|dx \right|  \leq 
 	& \displaystyle  C\ep \left( \ep^{1-\frac{1}{r}} + \ep^{1-\frac{1}{r}} \right) \|\varphi\|_{W^{1,r}_0(\Omega_\ep)} \leq \displaystyle C\ep^{2-\frac{1}{r}} \|\varphi\|_{W^{1,r}_0(\Omega_\ep)}.
 	\end{array}$$
Since ${2-\frac{1}{r}}= {1+\frac{1}{r'}}$, this yields
 	\begin{equation}
 		\left|  \ep(\eta_0-\eta_{\infty})	\int_{ \Omega_\ep}\left|(1+\lambda|\mathbb{D}[  {\bf u}_\ep]|^2)^{{r\over 2}-1}\mathbb{D}[ {\bf u}_\ep]:\mathbb{D}[   \varphi]\right|dx \right|  \nonumber\\
 	\leq  C\ep^{1+\frac{1}{r'}} \|\varphi\|_{W^{1,r}_0(\Omega_\ep)}. \label{step3_1}
 	\end{equation}
 	Similarly for the rest of the terms, by writing
 	\begin{align}
 		\left| \ep  \eta_{\infty}\int_{\Omega_\ep} \mathbb{D}[  {\bf u}_\ep]:\mathbb{D}[  \varphi]\, dx\right|
 		& \leq C\ep \|\mathbb{D}[  {\bf u}_\ep]\|_{L^2( \Omega_\ep)^{3\times 3}} \|\mathbb{D}[  \varphi]\|_{L^2( \Omega_\ep)^{3\times 3}}  \leq C \ep^{\frac{3}{2}} \ep^{\frac{1}{2}-\frac{1}{r}}\|\mathbb{D}[  \varphi]\|_{L^r( \Omega_\ep)^{3\times 3}}\nonumber \\
 		& \leq C\ep^{2-\frac{1}{r}}\|\mathbb{D}[  \varphi]\|_{L^r( \Omega_\ep)^{3\times 3}}\leq C\ep^{1+\frac{1}{r'}}\|\varphi\|_{W^{1,r}_0(\Omega_\ep)^3}, \label{step3_2}\\
 		\left| \int_{\Omega_\ep}f'\cdot \varphi' \,dx \right| & \leq  C \ep^{\frac{1}{r'}}\ep \|D \varphi\|_{L^r(\Omega_\ep)^{3\times 3}}  \leq C\ep^{1+\frac{1}{r'}} \|\varphi\|_{W^{1,r}_0(\Omega_\ep)^3}. \label{step3_3}
 	\end{align}
 Returning to expression~\eqref{Fep2}, we deduce from previous estimates that
\begin{equation}\label{estimDpr}\|\nabla p_\ep\|_{W^{-1,r'}(\Omega_\ep)^3}\leq C\ep^{{1\over r'}+1},
\end{equation}
which is (\ref{estim_nabla_p}) for $q=r$ and $q'=r'$. This finishes the proof  the case $r>2$.

Finally, the estimates of $\widetilde p^\ep_1$ given in (\ref{esti_P}) and (\ref{esti_P_r}) follows from the estimates of $p^\ep_1$ by using the change of variables (\ref{dilatacion}).

\end{proof}

\section{Adaptation of the unfolding method}\label{sec:unfolding}
 The change of variables (\ref{dilatacion}) does not capture the oscillations of the domain $\widetilde\Omega_\ep$. To capture them,  we use an adaptation of the unfolding method (see \cite{Ciora, Ciora2} for more details) introduced to this context in \cite{Anguiano_SG}. Thus, given $(\widetilde{\varphi}_{\varepsilon},  \zeta_\ep, \widetilde \psi_\varepsilon) \in L^q(\widetilde \Omega_\ep)^3\times L^{q'}(\omega)\times L^{q'}(\widetilde \Omega_\ep)$, $1\leq q<+\infty$ and $1/q+1/q'=1$, we define $(\widehat{\varphi}_{\varepsilon}, \widehat \zeta_\ep, \widehat \psi_\varepsilon)\in L^q(\omega\times Z)^3\times L^{q'}(\omega\times Z')\times \times L^{q'}(\omega\times Z)$ by
\begin{equation}\label{phihat}
\begin{array}{l}
\displaystyle\widehat{\varphi}_{\varepsilon}(x^{\prime},z)=\widetilde{\varphi}_{\varepsilon}\left( {\varepsilon^\ell}\kappa\left(\frac{x^{\prime}}{{\varepsilon^\ell}} \right)+{\varepsilon^\ell}z^{\prime},z_3 \right),\quad \hbox{a.e. }(x',z)\in \omega\times Z,\\
\displaystyle\widehat{\zeta}_{\varepsilon}(x^{\prime},z')=\zeta_{\varepsilon}\left( {\varepsilon^\ell}\kappa\left(\frac{x^{\prime}}{{\varepsilon^\ell}} \right)+{\varepsilon^\ell}z^{\prime} \right),\quad \hbox{a.e. }(x',z')\in \omega\times Z',\\
\noame
\displaystyle 
\widehat{\psi}_{\varepsilon}(x^{\prime},z)=\widetilde{\psi}_{\varepsilon}\left( {\varepsilon^\ell}\kappa\left(\frac{x^{\prime}}{{\varepsilon^\ell}} \right)+{\varepsilon^\ell}z^{\prime},z_3 \right),\quad \hbox{ a.e. }(x^{\prime},z)\in \omega\times Z,
\end{array}\end{equation}
assuming $\widetilde \varphi_\varepsilon$, $\zeta_\ep$ and $\widetilde \psi_\varepsilon$ are extended by zero outside $\omega$, where the function $\kappa:\mathbb{R}^2\to \mathbb{Z}^2$ is defined by 
$$\kappa(x')=k'\Longleftrightarrow x'\in Z'_{k',1},\quad\forall\,k'\in\mathbb{Z}^2.$$

\begin{remark}\label{remarkCV}We make the following comments:
\begin{itemize}
\item The function $\kappa$ is well defined up to a set of zero measure in $\mathbb{R}^2$ (the set $\cup_{k'\in \mathbb{Z}^2}\partial Z'_{k',1}$). Moreover, for every $\varepsilon>0$, we have
$$\kappa\left({x'\over \varepsilon^\ell}\right)=k'\Longleftrightarrow x'\in Z'_{k',\varepsilon^\ell}.$$

\item For $k^{\prime}\in \mathcal{T}_{\varepsilon}$, the restriction of $(\widehat{\varphi}_{\varepsilon},   \widehat \psi_\varepsilon)$  to $Z^{\prime}_{k^{\prime},{\varepsilon^\ell}}\times Z$ does not depend on $x^{\prime}$, whereas as a function of $z$ it is obtained from $(\widetilde\varphi_{\varepsilon},  \widetilde\psi_{\varepsilon})$ by using the change of variables 
\begin{equation}\label{CVunfolding}
\displaystyle z^{\prime}=\frac{x^{\prime}- {\varepsilon^\ell}k^{\prime}}{{\varepsilon^\ell}},\end{equation}
which transforms $Z_{k^{\prime}, {\varepsilon^\ell}}$ into $Z$. Analogously, the restriction of $\widehat \zeta_\ep$ to $Z'_{k',\ep^\ell}\times Z'$ does not depend on $x'$, while as function of $z'$ it is obtained from $\zeta_\ep$ by using the previous change of variables.
\end{itemize}
\end{remark}
Following the  proof of \cite[Lemma 4.9]{Anguiano_SG}, the following estimates relate  $(\widehat \varphi_\varepsilon,\widehat \zeta_\ep, \widehat \psi_\varepsilon)$ to $(\widetilde \varphi_\varepsilon, \zeta_\ep, \widetilde \psi_\varepsilon)$.
\begin{lemma}\label{estimates_relation}
We have the following estimates:
\begin{itemize}
\item  For every $\widetilde\varphi_\varepsilon \in L^q(\widetilde\Omega_\varepsilon)^3$, $1\leq q<+\infty$,
$$\|\widehat \varphi_\varepsilon\|_{L^q(\omega\times Z)^3}=  \|\widetilde \varphi_\varepsilon\|_{L^q(\widetilde\Omega_\varepsilon)^3},$$
where $\widehat \varphi_\varepsilon$ is given by (\ref{phihat})$_1$. Similarly, for every $\zeta_\ep\in L^{q'}(\omega)$ and $\widetilde\psi_\ep \in L^{q'}(\Omega)$, the functions $\widehat \zeta_\ep$ and $\widehat \psi_\varepsilon$, given by (\ref{phihat})$_{2,3}$ respectively satisfy
$$\|\widehat \zeta_\varepsilon\|_{L^{q'}(\omega\times Z')} = \|\zeta_\varepsilon\|_{L^{q'}(\omega)},\quad \|\widehat \psi_\varepsilon\|_{L^{q'}(\omega\times Z)} = \|\widetilde \psi_\varepsilon\|_{L^q(\widetilde \Omega_\ep)}.$$
\item  For every $\widetilde \varphi_\ep\in W^{1,q}(\widetilde\Omega_\varepsilon)^3$, $1\leq q<+\infty$, the function $\widehat \varphi_\varepsilon$ given by (\ref{phihat})$_1$ belongs to $L^q(\omega;W^{1,q}(Z)^3)$, and
$$\|D_{z'} \widehat \varphi_\varepsilon\|_{L^q(\omega\times Z)^{3\times 2}} = \varepsilon^\ell  \|D_{x'}\widetilde \varphi_\varepsilon\|_{L^q(\widetilde\Omega_\varepsilon)^{3\times 2}},\quad \|\partial_{z_3} \widehat \varphi_\varepsilon\|_{L^q(\omega\times Z)^{3 }} =  \|\partial_{z_3}\widetilde \varphi_\varepsilon\|_{L^q(\widetilde\Omega_\varepsilon)^{3}},$$
$$ \|\mathbb{D}_{z'}[\widehat \varphi_\varepsilon]\|_{L^q(\omega\times Z)^{3\times 2}} = \varepsilon^\ell  \|\mathbb{D}_{x'}[\widetilde \varphi_\varepsilon]\|_{L^q(\widetilde\Omega_\varepsilon)^{3\times 2}},\quad \|\partial_{z_3}[\widehat \varphi_\varepsilon]\|_{L^q(\omega\times Z)^{3 }} = \|\partial_{z_3}[\widetilde \varphi_\varepsilon]\|_{L^q(\widetilde\Omega_\varepsilon)^{3}}.$$
\item For every $\zeta_\ep\in W^{1,q'}(\omega_\ep)^3$, $1\leq q'<+\infty$, the function $\widehat \zeta_\ep$ given by (\ref{phihat})$_{2}$ belongs to $L^{q'}(\omega;W^{1,q'}(Z')^3)$, and 
$$\|\nabla_{z'}\widehat \zeta_\ep\|_{L^{q'}(\omega\times Z')^2}=\ep^{\ell}\|\nabla_{x'}\zeta_\ep\|_{L^{q'}(\omega)^2}.$$
\end{itemize} 
\end{lemma}

\begin{definition}[Unfolded velocity and pressures] Let us define the unfolded velocity and pressures $(\widehat {\bf u}_\varepsilon, \widehat p^0_\varepsilon, \widehat p^1_\ep)$ from $(\widetilde {\bf u}_\varepsilon, p^\varepsilon_0, \widetilde p^\ep_1)$ depending on the value of $r$:
\begin{itemize}
\item If $1<r< 2$,  from $(\widetilde {\bf u}_\varepsilon, p^\varepsilon_0, \widetilde p^1_\ep)\in H^1_0(\widetilde \Omega_\ep)^3\times H^1(\omega)\times L^2(\widetilde \Omega_\ep)$, we define $(\widehat {\bf u}_\varepsilon, \widehat p^\ep_0, \widehat p^1_\varepsilon)$ by (\ref{phihat}) with $\widetilde \varphi_\varepsilon=\widetilde {\bf u}_\varepsilon$, $\zeta_\ep=p^\ep_0$,  $\widetilde \psi_\varepsilon=\widetilde p\varepsilon_1$ and $q=q'=2$.\\

\item If $r>2$, from $(\widetilde {\bf u}_\varepsilon, p^\ep_0, \widetilde p^\varepsilon_1)\in W^{1,r}_0(\widetilde \Omega_\ep)^3\times W^{1,r'}(\omega)\times L^{r'}_0(\Omega)$, we define $(\widehat {\bf u}_\varepsilon, \widehat p_\ep^0, \widehat p_\varepsilon^1)$ by (\ref{phihat}) with $\widetilde \varphi_\varepsilon=\widetilde {\bf u}_\varepsilon$, $\zeta_\ep=p^\ep_0$, $\widetilde \psi_\varepsilon=\widetilde P_\varepsilon$, $q=r$ and $q'=r'$.
\end{itemize}
\end{definition}
Now, combining  estimates on the extended velocity (\ref{estim_sol_dil1})-(\ref{estim_sol_dil1_r}) and pressure (\ref{esti_P})-(\ref{esti_P_r}) with  Lemma~\ref{estimates_relation}, we deduce the following estimates on $({\bf \widehat u}_\varepsilon,\widehat p^0_\varepsilon, \widehat p^1_\ep)$.

\begin{lemma}\label{estimates_hat} The unfolded velocity/pressures triplet $(\widehat {\bf u}_\varepsilon,\widehat p^0_\ep, \widehat p^1_\varepsilon)$ satisfies the following estimates, depending on the value of $\ell$:
 \begin{itemize}
 \item If $1<r< 2$,  there exists a constant $C>0$ independent of $\varepsilon$, such that, 
  \begin{eqnarray}\medskip
 \|\widehat {\bf u}_\varepsilon\|_{L^2(\omega\times Z)^3}\leq C\varepsilon,& 
 \|D_{z'}\widehat {\bf u}_\varepsilon\|_{L^2(\omega\times Z)^{3\times 2}}\leq C\varepsilon^\ell,& 
  \|\partial_{z_3}\widehat {\bf u}_\varepsilon\|_{L^2(\omega\times Z)^{3}}\leq  C\varepsilon,\label{estim_u_hat}\\
  \medskip
 \|\widehat p^0_\varepsilon\|_{L^2(\omega\times Z')}\leq C,& \|\nabla_{z'}\widehat p^0_\varepsilon\|_{L^2(\omega\times Z')}\leq C\ep^\ell,&  \|\widehat p^1_\varepsilon\|_{L^2(\omega\times Z')}\leq C\ep.\label{estim_P_hat}
    \end{eqnarray}
 \item If $r>2$, there exists a constant $C>0$ independent of $\varepsilon$, such that,
 \begin{eqnarray}
 \|\widehat {\bf u}_\varepsilon\|_{L^r(\omega\times Z)^3}\leq C\varepsilon,& \|D_{z'}\widehat {\bf u}_\varepsilon\|_{L^r(\omega\times Z)^{3\times 2}}\leq C\varepsilon^\ell,& 
  \|\partial_{z_3}\widehat {\bf u}_\varepsilon\|_{L^r(\omega\times Z)^{3}}\leq  C\varepsilon,\label{estim_u_hatr}\\
    \noame
   \|\widehat p^0_\varepsilon\|_{L^{r'}(\omega\times Z')}\leq C,& \|\nabla_{z'}\widehat p^0_\varepsilon\|_{L^{r'}(\omega\times Z')}\leq C\ep^\ell,&  \|\widehat p^1_\varepsilon\|_{L^{r'}(\omega\times Z')}\leq C\ep,\label{estim_P_hatr}
  \end{eqnarray}
  where $r'$ is the conjugate exponent of $r$.
 \end{itemize}

 \end{lemma}

\section{Convergences of velocity and pressures}\label{sec:convergences}
In this section, we analyze the asymptotic behavior of extended functions $(\widetilde {\bf u}_\ep, p^\ep_0, \widetilde p^\ep_1)$ and the corresponding unfolded functions $(\widehat {\bf u}_\ep, \widehat p^0_\ep, \widehat p^1_\ep)$, when $\ep$ tends to zero. As we will see, the asymptotic behavior of velocities depends on the roughness regime considered, i.e. {\it Reynolds roughness regime} ($0<\ell<1$) or {\it Stokes roughness regime} ($\ell=1$), which will lead to deriving different limit systems in the next section.
\\

We define the following sets for $1<q<+\infty$. Let $C^\infty_{\#}(Z)$ be the space of infinitely differentiable functions in $\mathbb{R}^3$ that are $Z'$-periodic. By $L^q_{\#}(Z)$,  we denote its completion in the norm $L^q(Z)$ and by $L^q_{0,\#}(Z)$  the space of functions in $L^q_{\#}(Z)$ with zero mean value. Moreover, we introduce
\begin{equation}\label{Vz3}
\begin{array}{c}
\displaystyle
V_{z_3}^q(\Theta)=\{\varphi\in L^q(\Theta)\ :\ \partial_{z_3}\varphi\in L^q(\Theta)\},\quad V_{z_3,\#}^q(\Theta)=\{\varphi\in L^q_\#(\Theta)\ :\ \partial_{z_3}\varphi\in L^q(\Theta)\}\\
\noame
\displaystyle V_{z_3,\#}^q(\omega\times Z)=\{\varphi\in L^q(\omega;L^q_\#(Z))\ :\ \partial_{z_3}\varphi\in L^q(\omega\times Z)\}.
\end{array}
\end{equation} 

\begin{lemma}[Convergences of velocities] \label{lem_conv_vel}Assume $q=\max\{2,r\}$ and $q'$ the conjugate exponent of $q$. The velocities $\widetilde {\bf u}_\ep$ and $\widehat {\bf u}_\ep$ satisfy the following convergence results depending on the value of $\ell$ ($0<\ell<1$ or $\ell=1$):
\begin{itemize}
\item In the {\it Reynolds roughness regime} ($0<\ell<1$), then there exist
\begin{itemize}
\item  $\widetilde {\bf u}\in V_{z_3}^q(\Omega)^3$ where $\widetilde {\bf u}=0$ on $\Gamma_0\cup \Gamma_1$ and $\widetilde u_3\equiv 0$, such that, up to a subsequence, 
\begin{equation}\label{conv_vel_tilde2}
\ep^{-1}\widetilde {\bf u}_\ep\rightharpoonup \widetilde {\bf u}=(\widetilde {\bf u}',0)\quad\hbox{in }V_{z_3}^q(\Omega)^3,
\end{equation}
\item  $\widehat {\bf u}\in  V_{z_3,\#}^q(\omega\times Z)^3$, where $\widehat {\bf u}=0$ on $\omega\times (\widehat \Gamma_0\cup \widehat  \Gamma_1)$ and $\widehat u_3\equiv 0$, such that, up to a subsequence, 
\begin{equation}\label{conv_vel_hat2}
\ep^{-1}\widehat {\bf u}_\ep\rightharpoonup \widehat {\bf u}=(\widehat {\bf u}',0)\quad\hbox{in }V_{z_3}^q(\omega\times Z)^3,
\end{equation}
\begin{equation}\label{div_conv_z}
{\rm div}_{z'}(\widehat {\bf u}')=0\quad\hbox{in }\omega\times Z.
\end{equation}
\end{itemize}

\item In the {\it Stokes roughness regime} ($\ell=1$), then there exist
\begin{itemize}
\item  $\widetilde {\bf u}\in V_{z_3}^q(\Omega)$ where $\widetilde {\bf u}=0$ on $\Gamma_0\cup \Gamma_1$ and $\widetilde u_3\equiv 0$, such that, up to a subsequence, 
\begin{equation} \label{convStokestilde}
\ep^{-1}\widetilde {\bf u}_\ep\rightharpoonup (\widetilde {\bf u}',0)\quad\hbox{in }V_{z_3}^q(\Omega)^3,
\end{equation}
\item  $\widehat {\bf u}\in  L^q(\omega;W^{1,q}_\#(Z))^3$, where $\widehat {\bf u}=0$ on $\omega\times (\widehat \Gamma_0\cup \widehat  \Gamma_1)$ such that, up to a subsequence, 
\begin{equation} \label{convhatStokes}
\ep^{-1}\widehat {\bf u}_\ep\rightharpoonup \widehat {\bf u}\quad\hbox{in }L^q(\omega;W^{1,q}(Z))^3,
\end{equation}
\begin{equation}\label{div_conv_z_Stokes}
{\rm div}_{z}(\widehat {\bf u})=0\quad\hbox{in }\omega\times Z.
\end{equation}
\end{itemize}
\end{itemize}
Moreover,  in every case, it holds that 
\begin{equation}\label{div_conv}
{\rm div}_{x'}\left(\int_0^{h_{\rm max}}\widetilde {\bf u}'(x',z_3)\,dz_3\right)=0\quad \hbox{in }\omega,\quad\left(\int_0^{h_{\rm max}}\widetilde {\bf u}'(x',z_3)\,dz_3\right)\cdot n=0\quad \hbox{on }\partial\omega,
\end{equation}
and the following relation holds  
\begin{equation}\label{relation}
\widetilde {\bf u}(x',z_3)=\int_{Z'}\widehat {\bf u}(x',z)\,dz',\quad \hbox{with}\quad \int_{Z'}\widehat u_3(x',z)\,dz'=0,
\end{equation}
and so, 
\begin{equation}\label{relation2}
\int_0^{h_{\rm max}}\widetilde {\bf u}(x',z_3)\,dz_3=\int_{Z}\widehat {\bf u}(x',z)\,dz,\quad \hbox{with}\quad \int_{Z}\widehat u_3(x',z)\,dz=0.
\end{equation}
\end{lemma}
 \begin{proof}
Assume $q=\max\{2,r\}$. Let us first prove the case $0<\ell<1$.  To prove convergences (\ref{conv_vel_tilde2}), from (\ref{estim_sol_dil1}) and (\ref{estim_sol_dil1_r}), we have the following estimates
$$\|\widetilde{\bf u}_\ep\|_{L^q(\Omega)^3}\leq C\varepsilon,\quad \|D_{x'}\widetilde{\bf u}_\ep\|_{L^q(\Omega)^{3\times 2}}\leq C,\quad \|\partial_{z_3}\widetilde{\bf u}_\ep\|_{L^q(\Omega)^3}\leq C\varepsilon.$$
According to these estimates, we deduce that there exists $\widetilde {\bf u}\in V_{z_3}^q(\Omega)^3$ such that 
\begin{equation}\label{conv_vel_proof_1}
\ep^{-1}\widetilde {\bf u}_\ep\rightharpoonup \widetilde{\bf u}\quad\hbox{in }V_{z_3}^q(\Omega)^3,
\end{equation}
and moreover, $ \partial_{x_i}\widetilde {\bf u}_\ep$ tends to zero, $i=1,2$. Also, the continuity of the trace applications from the space of $\widetilde {\bf u}$ such that $\|\widetilde {\bf u}\|_{L^q}$ and $\|\partial_{z_3}\widetilde {\bf u}\|_{L^q}$  are bounded to $L^q(\Gamma_0\cup \Gamma_1)$ implies $\widetilde{\bf u}=0$ on $\Gamma_0\cup \Gamma_1$. Next, from the divergence condition ${\rm div}_\ep(\widetilde{\bf u}_\ep)=0$ in $\Omega$ and convergence (\ref{conv_vel_proof_1}), we deduce that $\partial_{z_3}\widetilde u_3=0$, i.e. $\widetilde u_3$ is independent of $z_3$. This together with $\widetilde u_3=0$ on $\Gamma_0\cup \Gamma_1$ implies that $\widetilde u_3\equiv 0$, which finishes the proof of (\ref{conv_vel_tilde2}).

By considering $\varphi\in C_c^1(\omega)$ as test function in divergence equation ${\rm div}_\ep(\widetilde{\bf u}_\ep)=0$ in $\Omega$, after integrating by parts, we deduce
$$\int_{\Omega}\widetilde{\bf u}'(x',z_3)\nabla_{x'}\varphi(x')\,dx'dz_3=0.$$
Multiplying by $\ep^{-1}$ and passing to the limit when $\ep$ tends to zero, since $\varphi$ does not depend on $z_3$, we easily deduce (\ref{div_conv}).

On the other hand, from estimates (\ref{estim_u_hat}) and (\ref{estim_u_hatr}), we have 
$$\|\widehat {\bf u}_\varepsilon\|_{L^q(\omega\times Z)^3}\leq C\varepsilon,\quad 
 \|D_{z'}\widehat {\bf u}_\varepsilon\|_{L^q(\omega\times Z)^{3\times 2}}\leq C\varepsilon^\ell,\quad
  \|\partial_{z_3}\widehat {\bf u}_\varepsilon\|_{L^q(\omega\times Z)^{3}}\leq  C\varepsilon.$$
Since $\ep\ll \ep^\ell$, by proceeding as for $\widetilde{\bf  u}_\ep$, we deduce convergences 
$$\ep^{-1}\widehat {\bf u}_\ep\rightharpoonup \widehat {\bf u}\quad\hbox{in }V_{z_3}^q(\omega\times Z)^3,$$
and $\ep^{-\ell}\partial_{z_i}\widehat {\bf u}_\ep$ tends to zero, $i=1,2$. Also, it holds that $\widehat {\bf u}=0$ on $\omega\times (\widehat \Gamma_0\cup \widehat \Gamma_1)$. By applying the unfolding change of variables to the divergence equation, we have that 
\begin{equation}\label{div_hat}\ep^{-\ell}{\rm div}_{z'}(\widehat{\bf u}'_\ep)+\ep^{-1}\partial_{z_3}\widehat u_3=0\quad\hbox{in }\omega\times Z,
\end{equation}
and passing to the limit, we deduce $\partial_{z_3}\widehat u_3=0$, and so $\widehat u_3\equiv 0$. It would remain to prove the $Z'$-periodicity of $\widehat {\bf u}$ in $z'$. This can be obtain by proceeding as in
\cite[Lemma 5.4]{grau1}. This finishes the proof of convergences (\ref{conv_vel_hat2}).

Now, taking a test function $\ep^{\ell -1}\widehat \varphi \in C^1_c(\omega\times Z')$ in (\ref{div_hat}), after integrating by parts, we deduce 
$$\int_{\omega\times Z}\ep^{-1}\widehat {\bf u}'\cdot \nabla_{z'}\widehat \varphi\,dx'dz=0,$$
and passing to the limit as $\ep$ tends to zero, after integrating by parts, we deduce (\ref{div_conv_z}).
\\

The proof of relation (\ref{relation}) is similar to the step 3 of the proof of \cite[Lemma 5.4]{Anguiano_SG}, so we omit it.
\\

 Finally, let us   prove the case $\ell=1$.   The proof of convergence (\ref{convStokestilde}) is similar to the previous case, so we omit it. On the other hand, from estimates (\ref{estim_u_hat}) and (\ref{estim_u_hatr}) with $\ell=1$, we have 
$$\|\widehat {\bf u}_\varepsilon\|_{L^q(\omega\times Z)^3}\leq C\varepsilon,\quad 
 \|D_{z'}\widehat {\bf u}_\varepsilon\|_{L^q(\omega\times Z)^{3\times 2}}\leq C\varepsilon,\quad
  \|\partial_{z_3}\widehat {\bf u}_\varepsilon\|_{L^q(\omega\times Z)^{3}}\leq  C\varepsilon.$$
Since $\ep\ll \ep^\ell$, by proceeding as for $\widetilde{\bf  u}_\ep$, we deduce convergences 
$$\ep^{-1}\widehat {\bf u}_\ep\rightharpoonup \widehat {\bf u}\quad\hbox{in }L^q(\omega;W^{1,q}(Z))^3,$$
and also, it holds that $\widehat {\bf u}=0$ on $\omega\times (\widehat \Gamma_0\cup \widehat \Gamma_1)$. By applying the unfolding change of variables to the divergence equation, we have that 
\begin{equation}\label{div_hatStokes} \ep^{-1}{\rm div}_{z}(\widehat{\bf u}_\ep)=0\quad\hbox{in }\omega\times Z,
\end{equation}
and passing to the limit, we deduce divergence condition (\ref{div_conv_z_Stokes}). The $Z'$-periodicity of $\widehat {\bf u}$ in $z'$  and the rest of properties can be deduced as in the previous case, so we omit it.

 \end{proof}

\begin{lemma}[Convergences of pressures]\label{lem_conv_press} Assume $q=\max\{2,r\}$ and $q'$ the conjugate exponent of $q$. Then, there exist $\widetilde p\in L^{q'}_0(\omega)\cap W^{1,q'}(\omega)$, independent of $z$, $\widehat p_0\in L^{q'}(\omega;W^{1,q'}_\#(Z'))$ and $\widehat p_1\in L^{q'}(\omega;L^{q'}_\#(Z))$ such that
\begin{equation}\label{pressuresH}
p^\ep_0\rightharpoonup \widetilde p\quad\hbox{in }W^{1,q'}(\omega),\quad p^\ep_0\to \widetilde p\quad\hbox{in }L^{q'}(\omega),\quad\ep^{-\ell}\nabla_{z'}\widehat p^0_\ep \rightharpoonup \nabla_{x'}\widetilde p+\nabla_{z'}\widehat p_0\quad\hbox{in }L^{q'}(\omega;L^{q'}(Z')),
\end{equation}
\begin{equation}\label{pressuresp1}
\ep^{-1}\widehat p^1_\ep\rightharpoonup\widehat p_1\quad\hbox{in }L^{q'}(\omega;L^{q'}(Z)).
\end{equation}
\end{lemma}

\begin{proof}
From estimates for $p^0_\ep$ and $\widehat p_\ep^0$, and the classical compactness result for the unfolding method for a bounded sequence in $W^{1,q'}$, see for instance \cite{Cioran-book}, we obtain convergences (\ref{pressuresH}).

Estimates for $\widehat p^1_\ep$ given in Lemma \ref{estimates_hat} imply the existence of $\widehat p_1\in L^{q'}(\omega; L^{q'}_{\#}(Z))$ such that up to a subsequence, convergence (\ref{pressuresp1}) holds. 

Since $\widetilde p_\ep$ has mean value zero in $\widetilde\Omega_\ep$, from the decomposition of the pressure and the unfolding change of variables, we deduce
$$0=\int_{\omega\times Z}\widehat p_\ep^0\,dx'dz'+\int_{\omega\times Z}\widehat p^1_\ep\,dx'dz.$$
Considering the convergence of $\widehat p^0_\ep$ to $\widetilde p$, that $\widehat p^1_\ep$ tends to zero, and that $\widetilde p$ does not depend on $z'$, we obtain
$$\int_\omega \widetilde p\,dx'=0,$$
and then, $\widetilde p$ has mean value zero in $\omega$.
\end{proof}

\section{Limit models}\label{sec:limitmodel}

Using monotonicity arguments together with Minty's lemma (see for instance \cite{Tapiero2, EkelandTemam}), we first derive in Section \ref{sec:ineq} a variational inequality, depending on the regime ($0<\ell<1$ or $\ell=1$), that will be useful in the proof of the main results. Next, we derive in Section \ref{sec:Reynolds} the limit system in the {\it Reynolds roughness regime} and   in Section \ref{sec:Stokes} the limit system in the {\it Stokes roughness regime}. Finally, we give some comments concerning the {\it high-frequency regime} (i.e. $\ell>1$) in Section \ref{sec:high}.
 
\subsection{Unfolding variational inequality}\label{sec:ineq} According to Lemma \ref{lem_conv_vel}, we choose a test function ${\bf v}(x',z)\in \mathcal{D}(\omega; C^\infty_{\#}(Z)^3)$  and boundary values   ${\bf v}'(x',z)=0$ on $\omega \times (\widehat \Gamma_0\cup \widehat \Gamma_1)$ and ${\rm div}_{x'}(\int_Z {\bf v}')\,dz=0$ in $\omega$ and $(\int_{Z}{\bf v}'\,dz )\cdot n=0$ on $\partial\omega$.  Moreover, in the case $0<\ell<1$, we will consider $v_3\equiv 0$ and ${\rm div}_{z'}({\bf v}')=0$ in $\omega\times Z$, while in the case $\ell=1$, we will consider ${\rm div}_{z}({\bf v})=0$ in $\omega\times Z$.\\

Multiplying (\ref{system_1_dil}) by ${\bf v}(x',x'/\varepsilon^\ell,z_3)$, integrating by parts, and taking into account the extension of $\widetilde {\bf u}_\varepsilon$ and the decomposition of the pressure, we have
$$\begin{array}{l}
\displaystyle
\medskip
\varepsilon (\eta_0-\eta_\infty)\int_{ \Omega }(1+\lambda|\mathbb{D}_\varepsilon[\widetilde {\bf u}_\varepsilon]|^2)^{{r\over 2}-1}\mathbb{D}_\varepsilon[\widetilde {\bf u}_\varepsilon] :\left( \mathbb{D}_{x'}[ {\bf v}]+\ep^{-\ell}\mathbb{D}_{z'}[{\bf v}] + \varepsilon^{-1}\partial_{z_3}[ {\bf v}]\right)\,dx'dz_3\\
\medskip
\displaystyle
+\varepsilon \eta_\infty\int_{ \Omega}\mathbb{D}_\varepsilon[\widetilde {\bf u}_\varepsilon] :\left( \mathbb{D}_{x'}[ {\bf v}']+\ep^{-\ell}\mathbb{D}_{z'}[{\bf v}] + \varepsilon^{-1}\partial_{z_3}[ {\bf v}]\right)\,dx'dz_3\\
\medskip
\displaystyle
-\int_{\widetilde \Omega_\ep}(p^0_\ep+\widetilde p^1_\varepsilon) \left({\rm div}_{x'}({\bf v}')+\ep^{-\ell} {\rm div}_{z'}({\bf v}')+\ep^{-1}\partial_{z_3}v_3\right)\,dx'dz_3=\int_{\Omega} {\bf f}'\cdot {\bf v}'\,dx'dz_3.
\end{array}
$$
By the change of variables given in Remark \ref{remarkCV}, we obtain
\begin{equation}\label{form_var_hat}\begin{array}{l}
\displaystyle
\medskip
 \ep(\eta_0-\eta_\infty)\int_{ \omega\times Z }(1+\lambda |\varepsilon^{-\ell} \mathbb{D}_{z'}[\widehat {\bf u}_\varepsilon] +\ep^{-1}\partial_{z_3}[\widehat {\bf u}_\ep]|^2)^{{r\over 2}-1}\left(\varepsilon^{-\ell} \mathbb{D}_{z'}[\widehat {\bf u}_\varepsilon] +\ep^{-1}\partial_{z_3}[\widehat {\bf u}_\ep]\right):( \ep^{-\ell}\mathbb{D}_{z'}[{\bf v}] + \ep^{-1}\partial_{z_3}[ {\bf v}]])\,dx'dz\\
\medskip
\displaystyle
+ \ep \eta_\infty\int_{ \omega\times Z}\left(\varepsilon^{-\ell} \mathbb{D}_{z'}[\widehat {\bf u}_\varepsilon] +\ep^{-1}\partial_{z_3}[\widehat {\bf u}_\ep]\right):( \ep^{-\ell}\mathbb{D}_{z'}[{\bf v}] + \ep^{-1}\partial_{z_3}[ {\bf v}])\,dx'dz\\
\medskip
\displaystyle
-\int_{\omega\times Z}(\widehat p_\ep^0+\widehat p^1_\varepsilon)\,\left({\rm div}_{x'}({\bf v}')+\varepsilon^{-\ell}{\rm div}_{z'}({\bf v}')+\ep^{-1}\partial_{z_3}v_3\right)\,dx'dz
=\int_{\omega\times Z} {\bf f}'\cdot {\bf v}'\,dx'dz+O_\varepsilon,
\end{array}
\end{equation}
where $O_\varepsilon$ is a generic real sequence depending on $\varepsilon$, that can change from line to line, and devoted to tend to zero.

Now, let us define the functional $J_r$ by
$$
J_r({\bf v})={\eta_0-\eta_\infty\over r\lambda}\int_{\omega\times Z}(1+\lambda|\varepsilon^{1-\ell} \mathbb{D}_{z'}[  {\bf v}] + \partial_{z_3}[ {\bf v}]|^2)^{r\over 2}dx'dz+{\eta_\infty\over 2}\int_{\omega\times Z}|\varepsilon^{1-\ell} \mathbb{D}_{z'}[  {\bf v}] + \partial_{z_3}[ {\bf v}]|^2dx'dz.
$$
Observe that $J_r$ is convex and Gateaux differentiable on $L^q(\omega;W^{1,q}_\#(Z)^3)$ with  $q=\max\{2,r\}$, (which is an adaptation of \cite[Proposition 2.1 and Section 3]{Baranger}) and $A_r=J'_r$ is given by
$$\begin{array}{l}\displaystyle
(A_r({\bf w}),{\bf v})\\
\noame
=\displaystyle (\eta_0-\eta_\infty)\int_{\omega\times Z}(1+\lambda|\varepsilon^{1-\ell} \mathbb{D}_{z'}[  {\bf w}] + \partial_{z_3}[ {\bf w}]|^2)^{{r\over 2}-1}(\varepsilon^{1-\ell} \mathbb{D}_{z'}[  {\bf w}] + \partial_{z_3}[ {\bf w}]):(\varepsilon^{1-\ell} \mathbb{D}_{z'}[  {\bf v}] + \partial_{z_3}[ {\bf v}])dx'dz\\
\noame
\displaystyle+\eta_\infty\int_{\omega\times Z}(\varepsilon^{1-\ell} \mathbb{D}_{z'}[  {\bf w}] + \partial_{z_3}[ {\bf w}]):(\varepsilon^{1-\ell} \mathbb{D}_{z'}[  {\bf v}] + \partial_{z_3}[ {\bf v}])dx'dz.
\end{array}
$$
Applying \cite[Proposition 1.1., p.158]{Lions2}, $A_r$ is monotone, \emph{i.e.}
\begin{equation}\label{monotonicity}
(A_r({\bf w})-A_r({\bf v}),{\bf w}-{\bf v})\ge0,\quad \forall {\bf w},{\bf v}\in L^q(\omega;W^{1,q}_\#(Z)^3).
\end{equation}
On the other hand, for all $\varphi\in \mathcal{D}(\omega; C^\infty_{\#}(Z)^3)$ and boundary value $\varphi=0$ on $\omega\times (\widehat \Gamma_0\cup \widehat \Gamma_1)$, satisfying the divergence conditions ${\rm div}_{x'}(\int_Z \varphi')\,dz=0$ in $\omega$ and $(\int_{Z}\varphi'\,dz )\cdot n=0$ on $\partial\omega$,  we choose ${\bf v}_\varepsilon$ defined by
\begin{equation}\label{testv}
{\bf v}_\varepsilon=\varphi-\varepsilon^{-1}\widehat {\bf u}_\varepsilon,
\end{equation}
as a test function in (\ref{form_var_hat}). Moreover, according to previous assumptions, in the case $0<\ell<1$, we will consider $\varphi_3\equiv 0$ and ${\rm div}_{z'}(\varphi')=0$ in $\omega\times Z$, while in the case $\ell=1$, we will consider ${\rm div}_{z}(\varphi)=0$ in $\omega\times Z$.

Then we obtain
\begin{equation*}\begin{array}{l}
\displaystyle
\medskip
  (A_r(\varepsilon^{-1}\widehat {\bf u}_\varepsilon),{\bf v}_\varepsilon)-\int_{\omega\times Z}(\widehat p_\ep^0+\widehat p^1_\varepsilon)\,\left({\rm div}_{x'}({\bf v}'_\varepsilon)+\varepsilon^{-\ell}{\rm div}_{z'}({\bf v}'_\varepsilon)+\ep^{-1}\partial_{z_3}v_{3,\ep}\right)\,dx'dz=\int_{\omega\times Z} {\bf f}'\cdot {\bf v}'_\varepsilon\,dx'dz+O_\varepsilon,
\end{array}
\end{equation*}
which is equivalent to
\begin{equation*}\begin{array}{l}
\displaystyle
\medskip
 (A_r(\varphi)-A_r(\varepsilon^{-1}\widehat {\bf u}_\varepsilon),{\bf v}_\varepsilon)-  (A_r(\varphi),{\bf v}_\varepsilon)
 \\
 \medskip
 \displaystyle+\int_{\omega\times Z}(\widehat p_\ep^0+\widehat p^1_\varepsilon)\,\left({\rm div}_{x'}({\bf v}'_\varepsilon)+\varepsilon^{-\ell}{\rm div}_{z'}({\bf v}'_\varepsilon)+\ep^{-1}\partial_{z_3}v_{3,\ep}\right)\,dx'dz\
 \medskip
 \displaystyle=-\int_{\omega\times Z} {\bf f}'\cdot {\bf v}'_\varepsilon\,dx'dz+O_\varepsilon.
\end{array}
\end{equation*}
Due to (\ref{monotonicity}), we can deduce
\begin{equation}\label{lim_var_nuevo}\begin{array}{l}
\displaystyle
\medskip
  (A_r(\varphi),{\bf v}_\varepsilon)-\int_{\omega\times Z}(\widehat p_\ep^0+\widehat p^1_\varepsilon)\,\left({\rm div}_{x'}({\bf v}'_\varepsilon)+\varepsilon^{-\ell}{\rm div}_{z'}({\bf v}'_\varepsilon)+\ep^{-1}\partial_{z_3}v_{3,\ep}\right)\,dx'dz\ge\int_{\omega\times Z} {\bf f}'\cdot {\bf v}'_\varepsilon\,dx'dz+O_\varepsilon.
\end{array}
\end{equation}
Therefore, depending on the value of $\ell$, we have:
\begin{itemize}
\item In the {\it Reynolds roughness regime} ($0<\ell<1$), since $\varphi_3\equiv 0$,    ${\rm div}_{z'}(\varphi')=0$  and  $\varepsilon^{-\ell}{\rm div}_{z'}(\widehat {\bf u}_\varepsilon')+\ep^{-1}\partial_{z_3}\widehat u_3=0$ in $\omega\times Z$, we have that
$$\ep^{-\ell}{\rm div}_{z'}(\varphi'-\ep^{-1}\widehat {\bf u}_\ep')+\ep^{-1}\partial_{z_3}(\varphi_3-\ep^{-1}\widehat {u}_{3,\ep})=\ep^{-\ell}{\rm div}_{z'}(\varphi')+\ep^{-\ell}{\rm div}_{z'}(\ep^{-1} \widehat {\bf u}_\ep')+\partial_{z_3}(\widehat \ep^{-1} u_{3,\ep})=0,$$
and so,  (\ref{lim_var_nuevo}) reads  as follows 
\begin{equation}\label{v_ineq_carreau}\begin{array}{l}
\displaystyle
\medskip
 (\eta_0-\eta_\infty)\int_{\omega\times Z}(1+\lambda|\varepsilon^{1-\ell} \mathbb{D}_{z'}[  \varphi'] + \partial_{z_3}[ \varphi']|^2)^{{r\over 2}-1}(\varepsilon^{1-\ell} \mathbb{D}_{z'}[  \varphi'] + \partial_{z_3}[ \varphi']):(\varepsilon^{1-\ell} \mathbb{D}_{z'}[  {\bf v}_\ep'] + \partial_{z_3}[ {\bf v}_\ep'])dx'dz\\
\noame
\medskip
\displaystyle+  \eta_\infty\int_{\omega\times Z}(\varepsilon^{1-\ell} \mathbb{D}_{z'}[  \varphi'] + \partial_{z_3}[ \varphi']):(\varepsilon^{1-\ell} \mathbb{D}_{z'}[  {\bf v}_\ep'] + \partial_{z_3}[ {\bf v}_\ep'])dx'dz\\
\medskip
\displaystyle-\int_{\omega\times Z}(\widehat p_\ep^0+\widehat p^1_\varepsilon)\,{\rm div}_{x'}({\bf v}'_\varepsilon)\,dx'dz\ge\int_{\omega\times Z} {\bf f}'\cdot {\bf v}'_\varepsilon\,dx'dz+O_\varepsilon.
\end{array}
\end{equation}
\item In the {\it Stokes roughness regime} ($\ell=1$), since ${\rm div}_{z}(\varphi)=0$ and $\ep^{-1}{\rm div}_{z}(\widehat {\bf u}_\ep)=0$ in $\omega\times Z$ , we have that
$$\ep^{-1}{\rm div}_{z'}(\varphi'-\ep^{-1}\widehat {\bf u}_\ep')+\ep^{-1}\partial_{z_3}(\varphi_3-\ep^{-1}\widehat {u}_{3,\ep})=\ep^{-1}{\rm div}_{z}(\varphi)+\ep^{-2}{\rm div}_{z}(\widehat {\bf u}_\ep)=0,$$
and so,  (\ref{lim_var_nuevo}) reads  as follows 
\begin{equation}\label{v_ineq_carreauStokes}\begin{array}{l}
\displaystyle
\medskip
 (\eta_0-\eta_\infty)\int_{\omega\times Z}(1+\lambda|\mathbb{D}_{z}[  \varphi]|^2)^{{r\over 2}-1}\mathbb{D}_{z}[  \varphi]: \mathbb{D}_{z}[ {\bf v}_\ep] \, dx'dz+  \eta_\infty\int_{\omega\times Z}\mathbb{D}_{z}[  \varphi] : \mathbb{D}_{z}[  {\bf v}_\ep] \, dx'dz\\
\medskip
\displaystyle-\int_{\omega\times Z}(\widehat p_\ep^0+\widehat p^1_\varepsilon)\,{\rm div}_{x'}({\bf v}'_\varepsilon)\,dx'dz\ge\int_{\omega\times Z} {\bf f}'\cdot {\bf v}'_\varepsilon\,dx'dz+O_\varepsilon.
\end{array}
\end{equation}
\end{itemize}
In previous inequalities, $O_\varepsilon$ is devoted to tends to zero. Next, by using previous variational inequalities, we give the main results of the paper, concerning the limit models, depending on the regime $0<\ell<1$ or $\ell=1$.
\subsection{Reynolds roughness regime}\label{sec:Reynolds}
We give the main result in the {\it Reynolds roughness regime} ($0<\ell<1$). We observe that the limit problem is the same for pseudoplastic fluids ($1<r<2$) and dilatant fluids ($r>2$), only the spaces change.
\begin{theorem} Let $q=\max\{2,r\}$ and $q'$ the conjugate exponent of $q$. The pair $(\widehat u, \widetilde p)\in V_{z_3,\#}^q(\omega\times Z)^2\times (L^{q'}_0(\omega)\cap W^{1,q'}(\omega))$ given in Lemmas \ref{lem_conv_vel} and \ref{lem_conv_press} is the unique solution of the following reduced two pressure problem:
\begin{equation}\label{limit_model}
\left\{\begin{array}{rl}
\displaystyle
-{1\over 2}\partial_{z_3}\left(  \left((\eta_0-\eta_\infty)\left(1+{\lambda\over 2}|\partial_{z_3} \widehat {\bf u}'|^2\right)^{{r\over 2}-1}+\eta_\infty\right)\partial_{z_3}\widehat {\bf u}'\right)+\nabla_{z'}\widehat \pi={\bf f}'(x')-\nabla_{x'}\widetilde p(x') &\hbox{in }\omega\times Z,\\
\noame
\displaystyle
{\rm div}_{z'}(\widehat {\bf u}')=0&\hbox{in }\omega\times Z,\\
\noame
\displaystyle
{\rm div}_{x'}\left(\int_Z\widehat u'\,dz\right)=0&\hbox{in }\omega,\\
\noame
\displaystyle \left(\int_Z\widehat u'\,dz\right)\cdot n=0&\hbox{on }\partial\omega,\\
\noame
\displaystyle \widehat u'=0&\hbox{on }\omega\times (\widehat\Gamma_0\cup\widehat\Gamma_1),\\
\noame
\displaystyle \widehat \pi\in L^{q'}(\omega;L^{q'}_{0,\#}(Z')).
\end{array}\right.
\end{equation}
\end{theorem}

\begin{proof} We develop the proof for every case $1<r<+\infty$, $r\neq 2$. To do this, let $q=\max\{2,r\}$. Observe that conditions (\ref{limit_model})$_{2, 3, 4, 5}$ are consequence of Lemma \ref{lem_conv_vel}. Now, let us pass to the limit in the variational inequality (\ref{v_ineq_carreau}) by taking into account (\ref{testv}), i.e. ${\bf v}_\varepsilon=\varphi-\varepsilon^{-1}\widehat {\bf u}_\varepsilon$.
\begin{itemize}
\item First term of (\ref{v_ineq_carreau}). From convergences given in Lemma \ref{lem_conv_vel} and relation (\ref{RelationAlpha}), we have
$$\begin{array}{l}
\displaystyle(\eta_0-\eta_\infty)\int_{\omega\times Z}(1+\lambda|\varepsilon^{1-\ell} \mathbb{D}_{z'}[  \varphi'] + \partial_{z_3}[ \varphi']|^2)^{{r\over 2}-1}(\varepsilon^{1-\ell} \mathbb{D}_{z'}[  \varphi'] + \partial_{z_3}[ \varphi']):(\varepsilon^{1-\ell} \mathbb{D}_{z'}[  {\bf v}_\ep'] + \partial_{z_3}[ {\bf v}_\ep'])dx'dz\\
\noame
=\displaystyle(\eta_0-\eta_\infty)\int_{\omega\times Z}(1+\lambda|\partial_{z_3}[ \varphi']|^2)^{{r\over 2}-1}\partial_{z_3}[ \varphi']:\partial_{z_3}[ \varphi'-\widehat {\bf u}'])\,dx'dz+O_\ep.
\end{array}$$
\item Second term of (\ref{v_ineq_carreau}). From convergences given in Lemma \ref{lem_conv_vel} and relation (\ref{RelationAlpha}), we have
$$\begin{array}{l}\displaystyle 
\eta_\infty\int_{\omega\times Z}(\varepsilon^{1-\ell} \mathbb{D}_{z'}[  \varphi'] + \partial_{z_3}[ \varphi']):(\varepsilon^{1-\ell} \mathbb{D}_{z'}[  {\bf v}_\ep'] + \partial_{z_3}[ {\bf v}_\ep'])dx'dz=
\eta_\infty\int_{\omega\times Z} \partial_{z_3}[ \varphi']:\partial_{z_3}[ \varphi'-\widehat{\bf u}']\,dx'dz+O_\ep.
\end{array}$$
\item Third term of (\ref{v_ineq_carreau}). From strong convergence of the pressure $p^0_\ep$, weak convergence of $p^1_\ep$ and weak convergence of the velocity given in Lemmas \ref{lem_conv_vel} and \ref{lem_conv_press}, we have
$$\begin{array}{l}\displaystyle \int_{\omega\times Z}\widehat p^0_\varepsilon\, {\rm div}_{x'}({\bf v}'_\varepsilon)\,dx'dz=\int_{\omega\times Z}\widetilde p(x')\, {\rm div}_{x'}(\varphi'-\widehat {\bf u}')\,dx'dz+O_\ep,\quad \left|\int_{\omega\times Z}\widehat p^1_\varepsilon\, {\rm div}_{x'}({\bf v}'_\varepsilon)\,dx'dz\right|\leq C\ep\to 0,
\end{array}$$
and so
$$  \int_{\omega\times Z}(\widehat p^0_\varepsilon+\widehat p^1_\ep)\, {\rm div}_{x'}({\bf v}'_\varepsilon)\,dx'dz=\int_{\omega\times Z}\widetilde p(x')\, {\rm div}_{x'}(\varphi'-\widehat {\bf u}')\,dx'dz+O_\ep.$$
\item Last term of (\ref{v_ineq_carreau}). From convergences given in Lemmas \ref{lem_conv_vel}, we have
$$\int_{\omega\times Z} {\bf f}'\cdot {\bf v}'_\varepsilon\,dx'dz= \int_{\omega\times Z} {\bf f}'\cdot (\varphi'-\widehat{\bf u}')\,dx'dz+O_\ep.$$
\end{itemize}
Taking into account that $\widetilde p$ does not depend on $z_3$ and conditions ${\rm div}_{x'}(\int_Z \varphi')\,dz=0$ and ${\rm div}_{x'}(\int_Z\widehat{\bf u}'dz)=0$, we deduce 
$$\int_{\omega\times Z}\widetilde p(x')\, {\rm div}_{x'}(\varphi'-\widehat {\bf u}')\,dx'dz=\int_\omega \widetilde p(x'){\rm div}_{x'}\left(\int_Z(\varphi'-\widehat {\bf u}')\,dz\right)dx'=0,$$
and so, from previous convergences, we have the limit variational inequality
$$\begin{array}{l}\displaystyle
(\eta_0-\eta_\infty)\int_{\omega\times Z}(1+\lambda|\partial_{z_3}[ \varphi']|^2)^{{r\over 2}-1}\partial_{z_3}[ \varphi']:\partial_{z_3}[ \varphi'-\widehat {\bf u}'])\,dx'dz\\
\noame
\displaystyle 
+\eta_\infty\int_{\omega\times Z} \partial_{z_3}[ \varphi']:\partial_{z_3}[\varphi'-\widehat{\bf u}']\,dx'dz\geq \int_{\omega\times Z} {\bf f}'\cdot (\varphi'-\widehat{\bf u}')\,dx'dz.
\end{array}$$
By Minty's lemma, this is equivalent to
\begin{equation}\label{limit_var}\begin{array}{l}\displaystyle
{1\over 2}(\eta_0-\eta_\infty)\int_{\omega\times Z}(1+{\lambda\over 2}|\partial_{z_3} \widehat {\bf u}'|^2)^{{r\over 2}-1}\partial_{z_3}\widehat {\bf u}':\partial_{z_3}{\bf v}'\,dx'dz
+{1\over 2}\eta_\infty\int_{\omega\times Z} \partial_{z_3}\widehat {\bf u}':\partial_{z_3}{\bf v}'\,dx'dz = \int_{\omega\times Z} {\bf f}'\cdot {\bf v}'\,dx'dz,
\end{array}
\end{equation}
which, by density, holds for every ${\bf v}'\in \mathcal{V}$ with
$$\mathcal{V}=\left\{\begin{array}{l}
\displaystyle {\bf v}'\in V_{z_3,\#}^q(\omega\times Z')\ :\ {\rm div}_{z'}({\bf v}')=0\quad\hbox{in }\omega\times Z\\
\noame
\displaystyle
{\rm div}_{x'}\left(\int_Z{\bf v}'\,dz\right)=0\quad\hbox{in }\omega,\quad \left(\int_Z{\bf v}'\,dz\right)\cdot n=0\quad\hbox{on }\partial\omega
\end{array}
\right\}.
$$
We observe that in (\ref{limit_var}), we have used that $$\partial_{z_3}[\widehat {\bf u}']: \partial_{z_3} [{\bf v}']={1\over 2}\partial_{z_3}\widehat {\bf u}' \cdot \partial_{z_3} {\bf v}'\quad \hbox{and}\quad |\partial_{z_3}[\widehat {\bf u}']|^2={1\over 2}|\partial_{z_3}\widehat {\bf u}'|^2.$$ 
Finally, reasoning as in \cite[Lemma 1.5]{Allaire}, the orthogonal of $\mathcal{V}$ is made of gradients of the form $\nabla_{x'}\widetilde\pi(x')+\nabla_{z'}\widehat \pi(x',z')$, with $\widetilde \pi\in L^{q'}_0(\omega)$ and $\widehat\pi(x',z')\in L^{q'}(\omega;L^{q'}_{0,\#}(Z'))$. Thus, integrating by parts, the limit variational formulation (\ref{limit_var}) is equivalent to the two-pressures reduced Stokes system (\ref{limit_model}). It remains to prove that $\widetilde \pi$ coincides with pressure $\widetilde p$. Thus can be easily done passing to the limit similarly as above by considering the test function $\varphi$, which is divergence-free only in $z'$, and by identifying limits. According to \cite[Propositions 3.2 and 3.3]{Tapiero2}, it can be proved that (\ref{limit_model}) has a unique solution  $(\widehat u', \widetilde p)\in V_{z_3,\#}^q(\omega\times Z)^2\times (L^{q'}_0(\omega)\cap W^{1,q'}(\omega)$), so the entire sequence converges.

\end{proof}
Let $q=\max\{2,r\}$ and $q'$ the conjugate exponent of $q$. Now, we give an approximation of the model (\ref{limit_model}), where the macroscopic scale is
  decoupled from the microscopic one. To do this, we consider $(\widehat {\bf w}_{\delta'}(z), \widehat q_{\delta'}(z'))$, for every $\delta'\in\mathbb{R}^2$, satisfying the local problem given by
\begin{equation}\label{limit_model_local}
\left\{\begin{array}{rl}
\displaystyle
-{1\over 2}\partial_{z_3}\left(  \left((\eta_0-\eta_\infty)\left(1+{\lambda\over 2}|\partial_{z_3} \widehat {\bf w}_{\delta'}|^2\right)^{{r\over 2}-1}+\eta_\infty\right)\partial_{z_3}\widehat {\bf w}_{\delta'}\right)+\nabla_{z'}\widehat q_{\delta'}(z')=-\delta' &\hbox{in }  Z,\\
\noame
\displaystyle
{\rm div}_{z'}(\widehat {\bf w}_{\delta'})=0&\hbox{in }  Z,\\
\noame
\displaystyle \widehat {\bf w}_{\delta'}=0&\hbox{on } \widehat\Gamma_0\cup\widehat\Gamma_1,\\
\noame
\displaystyle \widehat q_{\delta'}\in L^{q'}_{0,\#}(Z).
\end{array}\right.
\end{equation}
This problem has a unique solution $(\widehat {\bf w}_{\delta'}(z), \widehat q_{\delta'}(z'))\in V_{z_3,\#}^q(Z)^2\times (L^{q'}_{0,\#}(Z)\cap W^{1,q'}_\#(Z))$ for every $\delta'\in\mathbb{R}^2$ (this can be proved by adapting the proofs of \cite[Propositions 3.2 and 3.3]{Tapiero2}).  In fact, the velocity $\widehat {\bf w}_{\delta'}$ is given by
\begin{equation}\label{expresionw}
\widehat {\bf w}_{\delta'}(z)=-2\int_{{h(z')\over 2}-z_3}^{h(z')\over 2}{\xi\over \psi(2|\delta'+\nabla_{z'}\widehat q_{\delta'}(z')||\xi|)}d\xi\left(\delta'+\nabla_{z'}\widehat q_{\delta'}(z')\right),
\end{equation}
where $\psi$ is the inverse function of 
\begin{equation}\label{psifun}\tau=\zeta\sqrt{{2\over \lambda}\left\{{\zeta-\eta_\infty\over\eta_0-\eta_\infty}\right\}^{2\over r-2}-1},\end{equation}
which has a unique solution noted by $\zeta=\psi(\tau)$ for $\tau\in\mathbb{R}^+$,  and the pressure $\widehat q_{\delta'}$, by using the divergence-free condition ${\rm div}_{z'}(\int_{0}^{h(z')}\widehat {\bf w}_{\delta'}\,dz_3)=0$ in $Z'$, satisfies the following limit law
\begin{equation}\label{expressionq}
{\rm div}_{z'}\left\{\left(\int_{{-h(z')\over 2}}^{h(z')\over 2}{\left({h(z')\over 2}+\xi\right)\xi\over \psi(2|\delta'+\nabla_{z'}\widehat q_{\delta'}(z')||\xi|)}d\xi\right)\left(\delta'+\nabla_{z'}\widehat q_{\delta'}(z')\right)
\right\}=0\quad\hbox{in }Z'.
\end{equation}
Finally, we  give the expression for the filtration velocity 
\begin{equation}\label{FiltrationV_def}\widetilde {\bf V}(x')=\int_0^{h_{\rm max}}\widetilde {\bf u}(x',z_3)\,dz_3,
\end{equation} and the equation for pressure $\widetilde p$  in terms of the local problem, which is the main  result of the paper for the Reynols roughness regime ($0<\ell<1$).
\begin{corollary}\label{expresiones_finales}
The filtration velocity is given by
\begin{equation}\label{FiltrationVelcovcity}
\widetilde {\bf V}'(x')=-2\int_{Z'}\left(\int_{{-h(z')\over 2}}^{h(z')\over 2}{\left({h(z')\over 2}+\xi\right)\xi\over \psi(2|\nabla_{x'}\widetilde p(x')-{\bf f}'(x')+\nabla_{z'}\widehat q(z')||\xi|)}d\xi\right)\left(\nabla_{x'}\widetilde p(x')-{\bf f}'(x')+\nabla_{z'}\widehat q(z')\right)dz',
\end{equation}
and $\widetilde V_3\equiv 0$ in $\omega$. Here $\widehat q(z')$ satisfies (\ref{expressionq}) with $\delta'=\nabla_{x'}\widetilde p(x')-{\bf f}'(x')$.

Moreover, the Reynolds problem for pressure $\widetilde p\in L^{q'}_0(\omega)\cap W^{1,q'}(\omega)$ is given by
\begin{equation}\label{Reynolds_equation}{\rm div}_{x'}\widetilde {\bf V}'(x')=0\quad\hbox{in }\omega,\quad \widetilde {\bf V}'(x')\cdot n=0\quad\hbox{on }\partial\omega.
\end{equation}

\end{corollary}
\begin{proof}
Using the idea from \cite{Bourgeat1} to decouple the homogenized problem of Carreau type (\ref{limit_model}), for every $\delta'\in\mathbb{R}^2$ we consider the function $\mathcal{U}:\mathbb{R}^2\to \mathbb{R}^2$ given by
\begin{equation}\label{permeabilityU}\mathcal{U}(\delta)=\int_Z \widehat {\bf w}_{\delta'}(z)\,dz,
\end{equation}
where $\widehat {\bf w}_{\delta'}(z)$ is the solution of (\ref{limit_model_local}). Thus, $(\widehat {\bf u}',\widehat \pi)$ takes the form
$$\widehat{\bf u}'(x',z)=\widehat {\bf w}_{\nabla_{x'}\widetilde p(x')-{\bf f}'(x')}(z),\quad \widehat \pi(x',z)=\widehat q_{\nabla_{x'}\widetilde p(x')-{\bf f}'(x')}(z)\quad \hbox{in }\omega\times Z,$$
and then, from the relation $\widetilde V'(x')=\int_0^{h_{\rm max}}\widetilde {\bf u}'(x',z_3)\,dz_3=\int_Z\widehat {\bf u}'(x',z)\,dz$ given in (\ref{relation2}), and taking into account that  $\widetilde V_3(x')=\int_0^{h_{\rm max}}\widetilde u_3(x',z_3)\,dz_3=\int_Z\widehat u_3\,dz=0$, we deduce the filtration velocity 
$$\widetilde V'(x')=\int_0^{h_{\rm max}}\widetilde u'(x',z_3)\,dz_3=\mathcal{U}(\nabla_{x'}\widetilde p(x')-{\bf f}'(x')),\quad \widetilde V_3\equiv 0,\quad \hbox{in }\omega.$$
Then, by using (\ref{permeabilityU}) together with the expression (\ref{expresionw}) with $\delta'=\nabla_{x'}\widetilde p(x')-{\bf f}'(x')$, and taking into account that  (see \cite[Proof of Proposition 3.3]{Tapiero2} for more details)
$$\int_0^{h(z')}\widehat {\bf w}_{\nabla_{x'}\widetilde p(x')-{\bf f}'(x')}(z)\,dz_3= -2\int_{{-h(z')\over 2}}^{h(z')\over 2}(\nabla_{x'}\widetilde p(x')-{\bf f}'(x')+\nabla_{z'}\widehat q(z')){\left({h(z')\over 2}+\xi\right)\xi\over \psi(2|\nabla_{x'}\widetilde p(x')-{\bf f}'(x')+\nabla_{z'}\widehat q(z')||\xi|)}d\xi,$$
we deduce the expression for the filtration velocity (\ref{FiltrationVelcovcity}).

Finally, from the conditions (\ref{limit_model})$_{3,4}$, we deduce the Reynolds problem for pressure (\ref{Reynolds_equation}).
 
\end{proof}

\subsection{Stokes roughness regime}\label{sec:Stokes}
We give the main result in the {\it Stokes roughness regime} ($\ell=1$). We observe that the limit problem is the same for pseudoplastic fluids ($1<r<2$) and dilatant fluids ($r>2$), only the functional spaces are different.
\begin{theorem} Let $q=\max\{2,q\}$ and $q'$ the conjugate exponent of $q$. The pair $(\widehat u, \widetilde p)\in L^q(\omega;W^{1,q}(Z))^3\times (L^{q'}_0(\omega)\cap W^{1,q'}(\omega))$ given in Lemmas \ref{lem_conv_vel} and \ref{lem_conv_press} is the unique solution of the following reduced two pressure problem:
\begin{equation}\label{limit_modelStokes}
\left\{\begin{array}{rl}
\displaystyle
- {\rm div}(\eta_r(\mathbb{D}_z[\widehat{\bf u}])\mathbb{D}_z[\widehat {\bf u}])+\nabla_{z'}\widehat \pi={\bf f}'(x')-\nabla_{x'}\widetilde p(x') &\hbox{in }\omega\times Z,\\
\noame
\displaystyle
{\rm div}_{z}(\widehat {\bf u})=0&\hbox{in }\omega\times Z,\\
\noame
\displaystyle
{\rm div}_{x'}\left(\int_Z\widehat {\bf u}'\,dz\right)=0&\hbox{in }\omega,\\
\noame
\displaystyle \left(\int_Z\widehat {\bf u}'\,dz\right)\cdot n=0&\hbox{on }\partial\omega,\\
\noame
\displaystyle \widehat u=0&\hbox{on }\omega\times (\widehat\Gamma_0\cup\widehat\Gamma_1),\\
\noame
\displaystyle \widehat \pi\in L^{q'}(\omega;L^{q'}_{0,\#}(Z)),
\end{array}\right.
\end{equation}
where $\eta_r$ is defined by (\ref{Carreaulaw}), i.e.
\begin{equation}\label{CarreaulawStokes}\eta_r(\mathbb{D}_z[\widehat{\bf u}])=(\eta_0-\eta_\infty)(1+\lambda|\mathbb{D}_z[\widehat{\bf u}]|^2)^{{r\over 2}-1}+\eta_\infty,\quad 1<r<+\infty,\ r\neq 2,\quad\eta_0>\eta_\infty>0,\quad \lambda>0.
\end{equation}
\end{theorem}

\begin{proof} We develop the proof for every case $1<r<+\infty$, $r\neq 2$. To do this, let $q=\max\{2,r\}$. Observe that conditions (\ref{limit_modelStokes})$_{2, 3, 4, 5}$ are consequence of Lemma \ref{lem_conv_vel}. Now, let us pass to the limit in the variational inequality (\ref{v_ineq_carreauStokes}) by taking into account (\ref{testv}), i.e. ${\bf v}_\varepsilon=\varphi-\varepsilon^{-1}\widehat {\bf u}_\varepsilon$.
\begin{itemize}
\item First and second terms of (\ref{v_ineq_carreauStokes}). From convergences given in Lemma \ref{lem_conv_vel} and relation (\ref{RelationAlpha}), we have
$$\begin{array}{l}
\displaystyle (\eta_0-\eta_\infty)\int_{\omega\times Z}(1+\lambda|\mathbb{D}_{z}[  \varphi]|^2)^{{r\over 2}-1}\mathbb{D}_{z}[  \varphi]: \mathbb{D}_{z}[ {\bf v}_\ep] \, dx'dz+  \eta_\infty\int_{\omega\times Z}\mathbb{D}_{z}[  \varphi] : \mathbb{D}_{z}[  {\bf v}_\ep] \, dx'dz\\
\noame
=\displaystyle(\eta_0-\eta_\infty)\int_{\omega\times Z}(1+\lambda|\mathbb{D}_{z}[  \varphi]|^2)^{{r\over 2}-1}\mathbb{D}_{z}[  \varphi]: \mathbb{D}_{z}[ \varphi-\widehat{\bf u}] \, dx'dz+  \eta_\infty\int_{\omega\times Z}\mathbb{D}_{z}[  \varphi] : \mathbb{D}_{z}[  \varphi-\widehat{\bf u}] \, dx'dz+O_\ep.
\end{array}$$
\item Third term of (\ref{v_ineq_carreauStokes}). From strong convergence of the pressure $p^0_\ep$, weak convergence of $p^1_\ep$ and weak convergence of the velocity given in Lemmas \ref{lem_conv_vel} and \ref{lem_conv_press}, proceeding similarly to the previous regume, we deduce
$$  \int_{\omega\times Z}(\widehat p^0_\varepsilon+\widehat p^1_\ep)\, {\rm div}_{x'}({\bf v}'_\varepsilon)\,dx'dz=\int_{\omega\times Z}\widetilde p(x')\, {\rm div}_{x'}(\varphi'-\widehat {\bf u}')\,dx'dz+O_\ep.$$
\item Last term of (\ref{v_ineq_carreauStokes}). From convergences given in Lemmas \ref{lem_conv_vel}, we have
$$\int_{\omega\times Z} {\bf f}'\cdot {\bf v}'_\varepsilon\,dx'dz= \int_{\omega\times Z} {\bf f}'\cdot (\varphi'-\widehat{\bf u}')\,dx'dz+O_\ep.$$
\end{itemize}
Taking into account that $\widetilde p$ does not depend on $z_3$ and conditions ${\rm div}_{x'}(\int_Z \varphi')\,dz=0$ and ${\rm div}_{x'}(\int_Z\widehat{\bf u}'dz)=0$, we deduce 
$$\int_{\omega\times Z}\widetilde p(x')\, {\rm div}_{x'}(\varphi'-\widehat {\bf u}')\,dx'dz=\int_\omega \widetilde p(x'){\rm div}_{x'}\left(\int_Z(\varphi'-\widehat {\bf u}')\,dz\right)dx'=0,$$
and so, from previous convergences, we have the limit variational inequality
$$\begin{array}{l}\displaystyle
(\eta_0-\eta_\infty)\int_{\omega\times Z}(1+\lambda|\mathbb{D}_z[\varphi]|^2)^{{r\over 2}-1}\mathbb{D}_z[\varphi]:\mathbb{D}_z[ \varphi-\widehat {\bf u}])\,dx'dz\\
\noame
\displaystyle 
+\eta_\infty\int_{\omega\times Z} \mathbb{D}_z[\varphi]:\mathbb{D}_z[\varphi-\widehat{\bf u}]\,dx'dz\geq \int_{\omega\times Z} {\bf f}'\cdot (\varphi'-\widehat{\bf u}')\,dx'dz.
\end{array}$$
By Minty's lemma is equivalent to
\begin{equation}\label{limit_varStokes}\begin{array}{l}\displaystyle
(\eta_0-\eta_\infty)\int_{\omega\times Z}(1+{\lambda}|\mathbb{D}_z[\widehat {\bf u}]|^2)^{{r\over 2}-1}\mathbb{D}_z[\widehat {\bf u}]:\mathbb{D}_z[{\bf v}]\,dx'dz
+\eta_\infty\int_{\omega\times Z} \mathbb{D}_z[\widehat {\bf u}]:\mathbb{D}_z[{\bf v}]\,dx'dz = \int_{\omega\times Z} {\bf f}'\cdot {\bf v}'\,dx'dz,
\end{array}
\end{equation}
which, by density, holds for every ${\bf v}'\in \mathcal{V}$ with
$$\mathcal{V}=\left\{\begin{array}{l}
\displaystyle {\bf v}'\in L^q(\omega;W^{1,q}(Z)^3)\ :\ {\rm div}_{z}({\bf v})=0\quad\hbox{in }\omega\times Z\\
\noame
\displaystyle
{\rm div}_{x'}\left(\int_Z{\bf v}'\,dz\right)=0\quad\hbox{in }\omega,\quad \left(\int_Z{\bf v}'\,dz\right)\cdot n=0\quad\hbox{on }\partial\omega
\end{array}
\right\}.
$$
Finally, reasoning as in \cite[Lemma 1.5]{Allaire}, the orthogonal of $\mathcal{V}$ is made of gradients of the form $\nabla_{x'}\widetilde\pi(x')+\nabla_{z}\widehat \pi(x',yz)$, with $\widetilde \pi\in L^{q'}_0(\omega)$ and $\widehat\pi(x',z)\in L^{q'}(\omega;L^{q'}_{0,\#}(Z))$. Thus, integrating by parts, the limit variational formulation (\ref{limit_var}) is equivalent to the two-pressures   Stokes system (\ref{limit_modelStokes}). It remains to prove that $\widetilde \pi$ coincides with pressure $\widetilde p$. Thus can be easily done passing to the limit similarly as above by considering the test function $\varphi$, which is divergence-free only in $z'$, and by identifying limits. According to \cite[Propositions 3.2 and 3.3]{Tapiero2}, it can be proved that (\ref{limit_model}) has a unique solution  $(\widehat u', \widetilde p)\in L^q(\omega;W^{1,q}(Z)^3)\times (L^{q'}_0(\omega)\cap W^{1,q'}(\omega)$), so the entire sequence converges.

\end{proof}

Finally, we  give the expression for the filtration velocity ${\bf V}$ defined by (\ref{FiltrationV_def}) and the equation for pressure $\widetilde p$  in terms of the local problem, which is the main  result of the paper for the Stokes roughness regime ($\ell=1$).
\begin{corollary}\label{expresiones_finales_Stokes}
The filtration velocity is given by
\begin{equation}\label{FiltrationVelcovcityStokes}
\widetilde {\bf V}'(x')=\mathcal{U}({\bf f}'(x')-\nabla_{x'}\widetilde p(x')),\quad V_3(x')=0\quad\hbox{in }\omega,
\end{equation}
where the permeability function $\mathcal{U}:\mathbb{R}^2\to \mathbb{R}^2$ is defined by 
\begin{equation}\label{defU}\mathcal{U}(\xi')=\int_{Z}\widehat {\bf w}'_{\xi'}(z)\,dz,\quad \forall \xi'\in\mathbb{R}^2,
\end{equation}
with $\widehat {\bf w}_{\xi'}$, for every $\xi'\in\mathbb{R}^2$, the unique solution of the local Stoks system with non-linear viscosity following the Carreau law (\ref{Carreaulaw})
\begin{equation}\label{limit_modelStokes_local}
\left\{\begin{array}{rl}
\displaystyle
- {\rm div}(\eta_r(\mathbb{D}_z[\widehat {\bf w}_{\xi'}])\mathbb{D}_z[\widehat {\bf w}_{\xi'}])+\nabla_{z} \widehat \pi_{\xi'}=\xi' &\hbox{in } Z,\\
\noame
\displaystyle
{\rm div}_{z}(\widehat {\bf w}_{\xi'})=0&\hbox{in }Z,\\
\noame
\displaystyle \widehat {\bf w}_{\xi'}=0&\hbox{on }\widehat\Gamma_0\cup\widehat\Gamma_1,\\
\noame
\displaystyle \widehat {\bf w}_{\xi'}, \widehat \pi_{\xi'}\ Z'-periodic.
\end{array}\right.
\end{equation}
Moreover, the Reynolds problem for pressure $\widetilde p\in L^{q'}_0(\omega)\cap W^{1,q'}(\omega)$ is given by
\begin{equation}\label{Reynolds_equation_Stokes}{\rm div}_{x'}\widetilde {\bf V}'(x')=0\quad\hbox{in }\omega,\quad \widetilde {\bf V}'(x')\cdot n=0\quad\hbox{on }\partial\omega.
\end{equation}
\end{corollary}
\begin{proof}
To prove this result, we seek a global filtration velocity of the form 
\begin{equation}\label{filproof}
{\bf V}(x')=\mathcal{U}({\bf f}'(x')-\nabla_{x'}\widetilde p(x'))\quad\hbox{in }\omega,
\end{equation}
where $\mathcal{U}:\mathbb{R}^2\to \mathbb{R}^3$ is a permeability function, not necessarily linear, and ${\bf V}(x')=\int_0^{h_{\rm max}}\widetilde u(x',z_3)\,dz_3=\int_{Z}\widehat u(x',z)\,dz$ with ${\rm div}_{x'}{\bf V}'=0$ in $\omega$ and ${\bf V}'\cdot n=0$ on $\partial\omega$.

Using the idea from \cite{Bourgeat2} to decouple the homogenized problems of the Carreau type, for every $\xi'\in\mathbb{R}^2$ we consider the function $\mathcal{U}:\mathbb{R}^2\to \mathbb{R}^3$ given by
$$\mathcal{U}(\xi')=\int_{Z}\widehat w_{\xi'}(z)\,dz,$$
where $\widehat w_{\xi'}$ denotes the unique solution of the local Stokes problem given by (\ref{limit_modelStokes_local}), see \cite[Theorem 2]{Bourgeat1}. Thus, $(\widehat{\bf u},\widehat \pi)$ takes the form
$$\widehat{\bf u}(x',z)=\widehat w_{{\bf f}'(x')-\nabla_{x'}\widetilde p(x')},\quad\widehat\pi (x',z)=\widehat \pi_{{\bf f}'(x')-\nabla_{x'}\widetilde p(x')}\quad\hbox{in }\omega\times Z,$$
and then, from the relation ${\bf V}(x')=\int_{Z}\widehat {\bf u}(x',z)\,dz$ with $\int_{Z}\widehat u_3(x',z)\,dz=0$ given in (\ref{relation2}), we deduce the filtration velocity (\ref{filproof}), where $V_3=0$. Moreover, from(\ref{limit_modelStokes})$_{3, 4}$ together with (\ref{filproof}), we deduce (\ref{Reynolds_equation_Stokes}).

Since $V_3=0$, to simplify the notation, we redefine the definition of $\mathcal{U}$ by the expression given in (\ref{defU}) and then, we get $\mathcal{U}:\mathbb{R}^2\to \mathbb{R}^2$, which concludes the proof of (\ref{FiltrationVelcovcityStokes}). Finally, from \cite[Theorem 1]{Bourgeat2}, the macroscopic problem (\ref{FiltrationVelcovcityStokes})-(\ref{Reynolds_equation_Stokes}) has a unique solution $({\bf V}, \widetilde p)\in L^q(\omega)^3\times (L^q_0(\omega)\cap W^{1,q'}(\omega))$.

\end{proof}

\subsection{High-frequency roughness regime}\label{sec:high}
  In this section, we comment on the case of  {\it high-frequency roughness regime} ($\ell>1$), that is, the thickness of the domain $\varepsilon$ is greater than the period $\varepsilon^\ell$. We introduce the following notation:
\begin{itemize}
\item $\Omega_{+}^\varepsilon=\omega\times (\varepsilon h_{\rm min}, \varepsilon h(x'/\varepsilon^\ell))$ the oscillating part of the domain $\Omega_\varepsilon$,
\item $\widetilde \Omega_{+}^\varepsilon=\omega\times (h_{\rm min}, h(x'/\varepsilon^\ell))$ the oscillating part of the domain $\widetilde \Omega_\varepsilon$,
\item $\Omega_{+}=\omega\times (h_{\rm min}, h_{\rm max})$ the extension of the oscillating part,
\item $\widetilde {\bf u}_\varepsilon^+$ the restriction to $\widetilde \Omega_{+}^\varepsilon$ of the velocity. We denote by the same symbol the extension to $\Omega_+$, because the homogeneous boundary condition on the top boundary.

\end{itemize}

\noindent {\bf Estimates for velocity. } Because the following Poincar\'e's inequality in $\Omega_{+}^\varepsilon$ given in \cite[Lemma 4.3]{Anguiano_SG}:
$$\|w\|_{L^r(\Omega^\varepsilon_{+})}\leq C\varepsilon^\ell\|\nabla_{x'}w\|_{L^r(\Omega^\varepsilon_{+})^2},\quad \forall w\in W^{1,r}(\Omega_{+}^\varepsilon),\quad w_{|_{\widetilde \Gamma_1^\varepsilon}}=0,\quad 1\leq r<+\infty,$$
it can be proved that the following estimates hold:
$$
\|\widetilde {\bf u}_\varepsilon\|_{L^2(\widetilde \Omega^\varepsilon_{+})^3}\leq C\varepsilon^{2\ell-1} ,\quad \|D_\varepsilon \widetilde {\bf u}_\varepsilon\|_{L^2(\widetilde \Omega^\varepsilon_{+})^{3\times 3}}\leq C\varepsilon^{\ell-1},\quad \hbox{if}\quad1<r<+\infty,\ r\neq 2,
$$
and
$$\|\widetilde {\bf u}_\varepsilon\|_{L^r(\widetilde \Omega^\varepsilon_{+})^3}\leq C\varepsilon^{\ell+2{\ell-1\over r}}\quad \|D_\varepsilon \widetilde {\bf u}_\varepsilon\|_{L^r(\widetilde \Omega^\varepsilon_{+})^{3\times 3}}\leq C\varepsilon^{2{\ell-1\over r}}\quad \hbox{if} \quad r>2.$$
Observe that the same estimates also hold for the extensions to $\Omega_+$.\\

\noindent {\bf Convergences for velocity. } From the previous estimates for the extended $\widetilde {\bf u}_\varepsilon$ to $\Omega_+$, we deduce the following convergences
\begin{equation}\label{conv_u_high}
\varepsilon^{-\ell-2{\ell-1\over r}}\widetilde {\bf u}_\varepsilon\rightharpoonup 0\quad\hbox{ in }L^2(\Omega_+)^3\quad \hbox{if}\quad 1<r<2,\quad  
\varepsilon^{-\ell-2{\ell-1\over r}}\widetilde {\bf u}_\varepsilon\rightharpoonup 0\quad\hbox{ in }L^r(\Omega_+)^3 \quad \hbox{if}\quad r>2.
\end{equation}
The proofs respectively follow from the equalities
$$\varepsilon^{-1}\widetilde {\bf u}_\varepsilon=\varepsilon^{2(\ell-1)}\left(\varepsilon^{-2\ell+1} \widetilde {\bf u}_\varepsilon\right)\quad  \hbox{if}\quad 1<r<2,\quad \varepsilon^{-1}\widetilde {\bf u}_\varepsilon=\varepsilon^{(\ell-1)(1+{2\over r})}\left(\varepsilon^{-\ell-2{\ell-1\over r}} \widetilde {\bf u}_\varepsilon\right) \quad \hbox{if}\quad r>2,$$
and so, on the one hand,  since  $\varepsilon^{-2\ell+1}\widetilde {\bf u}_\varepsilon$ is bounded in $L^2(\Omega_+)$  and $\varepsilon^{2(\ell-1)} \to 0$, because $\ell>1$, we get the convergence for velocity in $L^2(\Omega_+)$, and    since $\varepsilon^{-\ell-2{\ell-1\over r}}  \widetilde {\bf u}_\varepsilon$ is bounded in $L^r(\Omega_+)$ and $\varepsilon^{(\ell-1)(1+{2\over r})} \to 0$, because $\ell>1$, we get the convergence for velocity in $L^r(\Omega_+)$. \\

\begin{figure}[h!]
\begin{center}
\includegraphics[width=15cm]{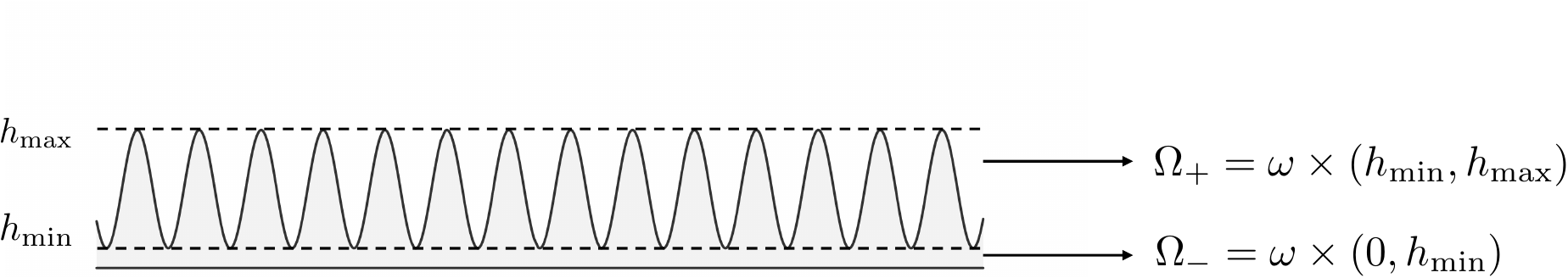}
\end{center}
\caption{View of  the subdomains of the high-frequency regime.}
\label{fig:omep_high}
\end{figure}

Therefore, the case $\ell>1$, due to the highly oscillating boundary, leads to the conclusion that
the velocity is zero in the roughness zone $\Omega_+$.  So, according to this,  at macroscopic level,  the problem reduces to the study of a Carreau fluid flow in the thin domain  $\Omega^\varepsilon_{-}=\omega\times (0,\varepsilon h_{\rm min})$, which after dilatation transform into a fixed set $\Omega_{-}=\omega\times (0, h_{\rm min})$, see Figure \ref{fig:omep_high}.\\

Thus, by applying the results of \cite{Tapiero2} to the Carreau fluid flow in the thin domain  $\Omega^\varepsilon_{-}=\omega\times (0,\varepsilon h_{\rm min})$, it follows that the limit model in $\Omega_{-}$  is given by 
\begin{equation}\label{limit_model_high}
\left\{\begin{array}{rl}
\displaystyle
-{1\over 2}\partial_{z_3}\left(  \left((\eta_0-\eta_\infty)\left(1+{\lambda\over 2}|\partial_{z_3} \widetilde{\bf u}'|^2\right)^{{r\over 2}-1}+\eta_\infty\right)\partial_{z_3}\widetilde {\bf u}'\right)+\nabla_{x'}\widetilde p(x') ={\bf f}'(x')&\hbox{in }\Omega_{-},\\
\noame
\displaystyle
{\rm div}_{x'}\left(\int_0^{h_{\rm min}}\widetilde{\bf u}'\,dz\right)=0&\hbox{in }\omega,\\
\noame
\displaystyle \left(\int_0^{h_{\rm min}}\widetilde{\bf u}'\,dz\right)\cdot n=0&\hbox{on }\partial\omega,\\
\noame
\displaystyle \widetilde {\bf u}'(x',0)=\widetilde {\bf u}'(x',h_{\rm min})=0,\quad \widetilde u_3\equiv 0.&
\end{array}\right.
\end{equation}
According to \cite[Proposition 3.3]{Tapiero2}, the  velocity $\widetilde {\bf u}'$ is given by
$$
\widetilde {\bf u}'(x',z_3)=2\int_{{h_{\rm min}\over 2}-z_3}^{h_{\rm min}\over 2}{\xi\over \psi(2|{\bf f}'(x')-\nabla_{x'}\widetilde p(x')||\xi|)}\,d\xi\left({\bf f}'(x')-\nabla_{x'}\widetilde p(x')\right),\quad (x',z_3)\in \Omega_{-},
$$
where $\psi$ is the inverse function of (\ref{psifun}), i.e. $\tau=\zeta\sqrt{{2\over \lambda}\left\{{\zeta-\eta_\infty\over\eta_0-\eta_\infty}\right\}^{2\over r-2}-1}$. Also, the filtration velocity $\widetilde {\bf V}(x')=\int_0^{h_{\rm min}}\widetilde u(x',z_3)\,dz_3$ is given by
$$
\widetilde {\bf V}'(x')=2\left({\bf f}'(x')-\nabla_{x'}\widetilde p(x')\right)\int_{-{h_{\rm min}\over 2}}^{h_{\rm min}\over 2}{\left({h_{\rm min}\over 2}+\xi\right)\xi\over \psi(2|{\bf f}'(x')-\nabla_{x'}\widetilde p(x')||\xi|)}\,d\xi,\quad \widetilde V_3\equiv 0\quad \hbox{in}\quad \omega.
$$
 Finally, the Reynolds problem for pressure $\widetilde p$ is given by
$${\rm div}_{x'}\widetilde {\bf V}'(x')=0\quad\hbox{in }\omega,\quad \widetilde {\bf V}'(x')\cdot n=0\quad\hbox{on }\partial\omega.
$$
 
\section{Conclusion}
This work provides a rigorous asymptotic analysis of the steady quasi-Newtonian Stokes flow, with viscosity described by the Carreau law, in a thin domain whose upper boundary exhibits small-scale oscillations. Through sharp a priori estimates, monotonicity arguments, and an adaptation of the unfolding method, we have derived effective two-dimensional nonlinear Reynolds-type equations that capture the influence of boundary roughness on the fluid flow.

The main results characterize how the interplay between the thickness parameter $\varepsilon$  and the wavelength of the boundary oscillations $\varepsilon^\ell$  governs the structure of the limit model, leading to distinct lubrication regimes:
\begin{itemize}
\item[--] {\it Reynolds Roughness Regime} ($0<\ell<1$): Here, the domain thickness is smaller than the characteristic wavelength of roughness. The limit model generalizes the classical Reynolds equation for Carreau fluids given in \cite{Tapiero2} by incorporating permeability functions reflecting the averaging effect of boundary oscillations over larger scales. We observe that the cell or local problems are of the same type of the limit equation obtained in \cite{Tapiero2}, but in the periodic setting. Thus, this results in a nonlinear Darcy-type formulation where the roughness alters the effective flow resistance in a non-trivial way, reflecting the non-Newtonian rheology.

\item[--] {\it Stokes Roughness Regime} ($\ell=1$): When the domain thickness and roughness wavelength scale proportionally, the rough boundary exerts more direct influence on the flow. The limit Reynolds equation includes microscopic cell problems of Stokes type that encode local boundary geometry and viscosity nonlinearities explicitly. This regime captures a more detailed coupling between geometry-induced microstructure and non-Newtonian viscosity effects, providing a refined description relevant for surfaces with moderate roughness.

\item[--] {\it High-frequency Roughness Regime} ($\ell>1$): Although not studied in full detail here, existing works suggest that in this regime the velocity tends to zero in the roughness zone due to highly oscillatory boundary conditions, effectively decoupling flow within the rough layer. This corresponds to a partial blocking effect of the roughness on the fluid motion.
\end{itemize}

Overall, the derived limit models demonstrate the significant role of boundary roughness in altering flow characteristics in non-Newtonian thin-film lubrication. In particular, the Carreau rheology leads to nonlinear Reynolds equations whose permeability depends intricately on both shear rates and roughness geometry. The mathematical framework developed herein, including rigorous compactness and homogenization techniques, sets a foundation for analyzing other complex fluids and boundary configurations.

From an applied perspective, these different regimes provide appropriate asymptotic descriptions for industrial and engineering lubrication problems involving non-Newtonian fluids with rough boundaries, allowing tailoring of models to the surface microstructure scale relative to lubricant thickness. The obtained results provide improved modeling tools for lubrication problems involving polymeric or shear-thinning fluids in engineering settings where surface microstructure cannot be neglected. Future work may extend these analyses to unsteady flows, temperature-dependent viscosities, or more irregular boundary roughness, broadening the applicability of nonlinear homogenized lubrication models.

 \section*{Acknowledgements}
Mar\'ia wants to dedicate this paper to her father, Julio, for all his love and for always accompanying her to the highest peaks. This work is based on preliminary developments presented in the preprint  \cite{Anguiano_SG_preprint} on the HAL repository. The authors would like to thank the reviewers for their helpful comments and suggestions that helped to improve the paper.

\section*{Author Contributions}
María Anguiano and Francisco J. Suárez-Grau have contributed equally to this work.

\section*{Funding} The authors declare that no funds, grants, or other support were received during the preparation of this manuscript.

 \section*{Data availability statement}
Data sharing not applicable to this article as no datasets were generated or analysed during the current study.

\section*{Competing Interests} The authors have no competing interests to declare that are relevant to the content of this article.

\section*{Conflict of interest} The authors state that there is no conflict of interest to declare.

\end{document}